\documentclass[
	english
]{scrartcl}

\usepackage[T1]{fontenc}
\usepackage[utf8]{inputenc}

\usepackage{amsmath,amsfonts,amsthm,amssymb}
\usepackage[colorlinks,citecolor=blue]{hyperref}
\usepackage{cancel}
\usepackage{relsize}
\usepackage{algorithmic}

\usepackage{changes}

\usepackage{enumitem}
\setlist[enumerate]{label=(\alph*)}

\usepackage[figure]{hypcap}
\usepackage{graphicx}
\graphicspath{{./img/}}
\usepackage{subcaption} 

\usepackage{preamble}

\usepackage{marginnote}

\usepackage[capitalise,nameinlink]{cleveref}
\allowdisplaybreaks

\numberwithin{equation}{section}

\crefname{assumption}{Assumption}{Assumptions}
\Crefname{ALC@unique}{Step}{Steps}

\newcommand\norm[1]{\left\Vert#1\right\Vert}
\newcommand\nnorm[1]{\Vert#1\Vert}

\newcommand\N{\mathbb{N}}
\newcommand\R{\mathbb{R}}

\newcommand\B{\mathbb{B}}

\newcommand\tto{\rightrightarrows}

\newcommand{\dist}{\operatorname{dist}}
\newcommand{\conv}{\operatorname{conv}}

\newcommand{\dom}{\operatorname{dom}}

\newcommand{\gph}{\operatorname{gph}}
\newcommand{\epi}{\operatorname{epi}}

\newcommand{\xb}{\bar x}

\renewcommand{\d}{\mathrm d}

\DeclareMathAlphabet{\mathpzc}{OT1}{pzc}{m}{it}
\newcommand\oo{\mathpzc{o}}

\newtheorem{theorem}{Theorem}[section]
\newtheorem{lemma}[theorem]{Lemma}
\newtheorem{proposition}[theorem]{Proposition}
\newtheorem{assumption}[theorem]{Assumption}
\newtheorem{algorithm}[theorem]{Algorithm}

\newtheorem{corollary}[theorem]{Corollary}
\newtheorem{remark}[theorem]{Remark}
\newtheorem{definition}[theorem]{Definition}
\newtheorem{example}[theorem]{Example}

\makeatletter
\long\def\@firstoffiveparen#1#2#3#4#5{\textup{\tagform@{#1}}}
\def\eqref@nolink#1{\textup{\tagform@{\ref*{#1}}}}
\def\eqref@link#1{%
\Hy@safe@activestrue
\expandafter\@setref\csname r@#1\endcsname\@firstoffiveparen{#1}%
\Hy@safe@activesfalse
}
\protected\def\eqref{\@ifstar\eqref@nolink\eqref@link}
\makeatother

\definecolor{mygreen}{rgb}{0.0,0.7,0.0}
\definecolor{mybrown}{rgb}{0.5,0.5,0.0}

\begin{document}

\title{%
	Isolated calmness of perturbation mappings in generalized nonlinear programming 
	and local superlinear convergence of Newton-type methods
	}%
\author{%
	Mat\'{u}\v{s} Benko%
	\footnote{%
		Johann Radon Institute for Computational and Applied Mathematics,
		4040 Linz,
		Austria
		\email{matus.benko@ricam.oeaw.ac.at},
		\orcid{0000-0003-3307-7939}
		}
	\and
	Patrick Mehlitz%
	\footnote{%
		Philipps-Universit\"at Marburg,
		Department of Mathematics and Computer Science,
		35032 Marburg,
		Germany,
		\email{mehlitz@uni-marburg.de},
		\orcid{0000-0002-9355-850X}%
		}
	}

\publishers{}
\maketitle

\begin{abstract}
	 In this paper, we characterize Lipschitzian properties of different
	 multiplier-free and multiplier-dependent perturbation mappings associated
	 with the stationarity system of a so-called generalized nonlinear program
	 popularized by Rockafellar. Special emphasis is put on the investigation
	 of the isolated calmness property at and around a point. 
	 The latter is decisive for the locally fast convergence of
	 the so-called semismooth* Newton-type method by Gfrerer and Outrata.
	 Our central result is the characterization of the isolated calmness at a point
	 of a multiplier-free perturbation mapping via a combination of an explicit condition
	 and a rather mild assumption, automatically satisfied e.g.\ for standard nonlinear programs. 
 	 Isolated calmness around a point is characterized analogously 
 	 by a combination of two stronger conditions.
	 These findings are then related to so-called criticality of
	 Lagrange multipliers, as introduced by Izmailov and extended to
	 generalized nonlinear programming by Mordukhovich and Sarabi.
	 We derive a new sufficient condition (a characterization for some problem classes) 
	 of nonexistence of critical multipliers,
	 which has been also used in the literature as an assumption to
	 guarantee local fast convergence of Newton-, SQP-, or multiplier-penalty-type methods.
	 The obtained insights about critical multipliers seem to complement 
	 the vast literature on the topic.
\end{abstract}


\begin{msc}	
	\mscLink{49J52}, \mscLink{49J53}, \mscLink{90C30}
\end{msc}

\section{Introduction}\label{sec:introduction}

The starting point for this paper has been a new interpretation of the condition 
defining so-called \emph{critical Lagrange multipliers},
which we noticed particularly after it was extended 
from the setting of standard nonlinear programs (NLPs), 
see \cite{Izmailov2005,IzmailovSolodov2012},
to more general settings by Mordukhovich and Sarabi, 
see \cite{MordukhovichSarabi2018,MordukhovichSarabi2019}.
Perhaps the most interesting feature of critical multipliers is that they cause 
a slow convergence of Newton-type methods,
see \cite{Izmailov2005,IzmailovKurennoySolodov2015,IzmailovSolodov2009,IzmailovSolodov2012,IzmailovSolodov2015}.
From the initial idea of interpreting critical multipliers, 
a broader motivation arose to better understand the regularity-type assumptions
known to be essential for a superlinear rate of convergence of Newton-type methods.
The novel \emph{semismooth*} Newton method (or rather this class of methods) 
from \cite{GfrererOutrata2021} seems to provide a particularly suitable framework for this task.

As the underlying optimization problem, we consider the model initiated in \cite{LevyPoliquinRockafellar2000}, 
where the problem is given as
\begin{equation}\label{eq:Pbar}\tag{P}
	\min\limits_x\quad f(x,0) \quad\text{s.t.}\quad x\in\R^n 
\end{equation}
for a proper, lower semicontinuous function $f\colon\R^n\times\R^m\to\overline{\R}$,
and it is interpreted as embedded in a parametrized family of problems
\begin{equation}\label{eq:P(u,v)}\tag{P$(v,u)$}
	\min\limits_x\quad f(x,u) - v ^\top x \quad\text{s.t.}\quad x\in\R^n
\end{equation}
involving parameters $v\in\R^n$ and $u\in\R^m$.
Particularly, \eqref{eq:Pbar} equals \hyperref[eq:P(u,v)]{(P$(0,0)$)}.
	Special emphasis is put on the situation where $f$ takes the composite form
\begin{equation}\label{eq:GNLP_setting}
f(x,u) := f_0(x) + g(F(x) + u)
\end{equation}
for twice continuously differentiable functions $f_0\colon\R^n\to\R$ and $F\colon\R^n\to\R^m$
as well as a proper, lower semicontinuous, convex function $g\colon\R^m\to\overline{\R}$. 
The corresponding composite optimization problem \eqref{eq:P(u,v)} 
has been referred to as a \emph{generalized} nonlinear programming problem
in \cite{BenkoRockafellar2023,Rockafellar2023,Rockafellar2023b} as it covers, exemplary,
standard nonlinear, 
nonlinear second-order cone,
and sparse optimization.
More details about the model can be found in \cref{sec:model}.

Essentially, the paper revolves around the set-valued mapping 
$M\colon\R^n\times\R^m\tto\R^n$ given by
\begin{equation*}
	M(v,u)
	:=
	\{x \,|\, v \in \partial_x f(x,u)\},
\end{equation*}
particularly around the properties of \emph{isolated calmness} and its new extension 
from \cite{GfrererOutrata2022,GfrererOutrata2023},
\emph{isolated calmness on a neighborhood}, of $M$ at/around 
a point $((0,0),\bar x)$ from its graph
for some given stationary point $\bar x\in\R^n$ of \eqref{eq:Pbar}.

Let us now summarize the main contributions of the paper, 
all valid under the composite structure \eqref{eq:GNLP_setting}.
\begin{itemize}
	\item We show that the isolated calmness of $M$ around $((0,0),\bar x)$ 
		can be used as an essential assumption to obtain local superlinear convergence 
		of a suitable variant of the semismooth* Newton method 
		to find a stationary point of \eqref{eq:Pbar}.\\
		This is actually just a small extension of the convergence theory 
		for the semismooth* Newton method from the papers
		\cite{GfrererOutrata2021,GfrererOutrata2022},
		but it provides an important justification 
		for focusing on the isolated calmness (on a neighborhood) of $M$.
	\item We fully characterize the isolated calmness of $M$ at $((0,0),\bar x)$ 
		via a combination of an explicit condition and a seemingly mild assumption.
	 	Similarly, a combination of certain robust versions of the explicit condition 
	 	and the mild assumption fully characterize 
	 	the isolated calmness of $M$ around $((0,0),\bar x)$.\\
		The mild assumption, necessary for the isolated calmness at the point,
		is the so-called \emph{inner calmness* in the fuzzy sense},
		and it is actually the key ingredient of the paper.
		It has been recently introduced in \cite{Benko2021} for its role 
		in the image rule for tangent cones, see \cite[Theorem~4.1]{Benko2021}.
		That calculus rule, in its extended version from \cite[Theorem 3.1]{BenkoMehlitz2022},
		is exactly what enables us to connect the isolated calmness of $M$ 
		with the explicit condition.
		For a broad class of problems where the epigraph of $g$ is convex polyhedral,
		covering NLPs but also many more, the inner calmness* assumption
		is automatically satisfied, 
		and the isolated calmness of $M$ is equivalent to the explicit condition.
		However, that explicit condition actually corresponds to nonexistence 
		of critical multipliers for $\bar x$.
		This brings us to the final and main contribution.
	\item We reveal the strong connection 
		between nonexistence of critical multipliers for $\bar x$ 
		and the isolated calmness of $M$ at $((0,0),\bar x)$.\\
		As we will show, 
		for composite problems modeled with a function $g$ whose epigraph is convex polyhedral
			and satisfying a qualification condition,
		nonexistence of critical multipliers is fully equivalent 
		with the isolated calmness of $M$.
			In general, the latter provides a sufficient condition for ruling out
			the existence of critical multipliers even without a qualification condition.
			This considerably improves the known result,
			established for a particular class of problems e.g.\ in \cite[Theorem~5.1]{MordukhovichSarabi2018},
			that so-called \emph{full stability} of $\bar x$ prevents critical multipliers.
		We believe that this connection 
			between critical multipliers and the isolated calmness of $M$
		provides an important, previously unnoticed insight into critical multipliers,
		and complements the vast literature on the topic.
\end{itemize}

Given the disruptive nature of critical multipliers,
ruling out their existence is the seemingly most preferable way to deal with them.
To the best of our knowledge, perhaps the only condition known to do so
in a rather general framework is the aforementioned full stability of $\bar x$.
In \cite[Theorem~5.1]{MordukhovichSarabi2018}, the authors prove this for convex, piecewise linear $g$.
It is known that the milder notion of \emph{tilt stability} is not enough;
in \cite[Corollary~6.1] {MordukhovichSarabi2018}, 
the authors rule out critical multipliers by tilt stability but only in combination with a certain 
nondegeneracy-type condition, which implies uniqueness of the multiplier.
The interesting degenerate case discussed in \cite[Remark~6.2]{MordukhovichSarabi2018} likely triggers full stability 
in the light of new insights from \cite[Theorem~5.35]{Mordukhovich2024}.
Our approach is significantly stronger as we prove that
the isolated calmness of $M$ rules out critical multipliers in the general composite setting \eqref{eq:GNLP_setting},
without a qualification condition,
and for many problem classes it is also necessary for the task.
Additionally, it is much milder than full stability 
(which actually implies the existence of a single-valued Lipschitz localization of $M$ around $((0,0),\bar x)$ 
as well as validity of a qualification condition).

Ruling out the existence of critical multipliers, however, is not always feasible
and a lot of effort has been put into understanding how to ``avoid'' them.
Central to such efforts is a finer, local approach, characterizing criticality of a single multiplier
in terms of an error bound condition, see e.g.\ \cite[Proposition~1]{IzmailovSolodov2012},
\cite[Theorem~4.1]{MordukhovichSarabi2018}, \cite[Theorem~5.6]{MordukhovichSarabi2019},
or \cite[Theorem~3.6, Proposition~3.8]{Sarabi2022}.
Interestingly, at the heart of these error bound characterizations
is a connection between critical multipliers and some isolated calmness assumption,
see e.g. \cite[Section~7]{MordukhovichSarabi2018}.
However, this isolated calmness assumption does not involve the mapping $M$,
but the mapping which assigns to parameters $(v,u)$ not only $x$, 
but also the corresponding multiplier $y$,
see the mapping $M_1$ in \cref{sec:model}.

In the literature, these error bounds were used to establish local fast convergence 
of Newton-, SQP-, or multiplier-penalty-type methods even in the case 
where critical multipliers exist.
In a forthcoming paper, we plan to adjust our approach to obtain this local kind of analysis
for critical multipliers for the general composite model \eqref{eq:CP},
and we will carve out some consequences of these findings for the local convergence of the
semismooth* Newton method from \cref{alg:SSN}.
For better context, let us refer the interested reader to various other Newton-type
methods for the numerical solutions of generalized and nonsmooth equations e.g.\ 
in the papers \cite{Josephy1979,KlatteKummer2018,Pang1990,QiSun1993,Robinson1994}
and the monographs \cite{IzmailovSolodov2014,KlatteKummer2002}.

The remainder of the paper is organized as follows.
In \cref{sec:preliminaries}, we comment on the notation in this manuscript 
and present some preliminary results. 
\cref{sec:basics} comprises comments about the fundamental notation we are using.
\cref{sec:tangents_and_normals} is dedicated to the introduction of several tangent and
normal cones, and we also present some calculus rules for the comparatively less popular 
limiting tangent and paratingent cone.
In \cref{sec:generalized_differentiation}, we discuss certain subdifferentials as well as
the second subderivative of proper lower semicontinuous functions.
The analysis of set-valued mappings is investigated in \cref{sec:set_valued_maps}.
We start by coining some standard terminology before introducing generalized derivatives
of set-valued mappings. We also recall several notions of regularity and Lipschitzness
of set-valued mappings and present some characterizations of these conditions in terms of 
the introduced derivatives. Special emphasis is laid on the fuzzy inner calmness*
property, and we also introduce a locally uniform version of it.
These two appear in our characterizations of the isolated calmness at/around a point of $M$.
We close \cref{sec:preliminaries} with a brief recapitulation of the semismooth*
Newton method from \cite{GfrererOutrata2021} 
for the numerical solution of generalized equations in \cref{sec:ss*}.
Particularly, we show the local superlinear convergence of this method 
under the strong metric subregularity (isolated calmness of the inverse) around a given solution of the generalized equation.
\cref{sec:perturbed_optimization_problems} of the paper is dedicated to the variational analysis
of certain perturbation mappings in nonlinear optimization.
In \cref{sec:model}, we provide details regarding the parametrized model problem \eqref{eq:P(u,v)}.
	Particularly, we discuss some beneficial consequences of the composite setting \eqref{eq:GNLP_setting}.
Afterwards, we introduce the perturbation mappings associated with the problem, and
some straightforward relations between those mappings are carved out.
\cref{sec:variational_analysis_perturbation_maps} is dedicated to
the variational analysis of the introduced perturbation mappings.
Particular focus is placed
on characterizations of the presence of the isolated calmness property (in a neighborhood)
of these mappings.
These findings are applied in \cref{sec:critical_multipliers} in order to state 
characterizations of nonexistence of critical multipliers.
We also introduce the novel concept of strong noncriticality and present similar
characterizations of it. Then we bridge the convergence results from \cref{sec:ss*},
which address the semismooth* Newton method, with the findings in 
\cite{Izmailov2005,IzmailovSolodov2012} regarding the Newton method 
in standard nonlinear optimization.
Some concluding comments close the paper in \cref{sec:conclusions}.

\section{Preliminaries}\label{sec:preliminaries}

In this paper, we mainly rely on standard notation as used in the textbooks 
\cite{DontchevRockafellar2014,Mordukhovich2018,RockafellarWets1998}.

\subsection{Fundamental notation}\label{sec:basics}

Throughout the paper, we use $\overline{\R}:=\R\cup\{-\infty,\infty\}$.
Furthermore, for $n\in\N$, $\R^n$ and $\R^n_-$ denote the sets of all real vectors with $n$ components
and $n$ nonpositive components, respectively.
We use $\mathtt e_1,\ldots,\mathtt e_n\in\R^n$ for the
$n$ canonical unit vectors of $\R^n$.
The set of all real matrices with $m\in\N$ rows and $n$ columns 
will be represented by $\R^{m\times n}$,
and for each $A\in\R^{m\times n}$, $A^\top\in\R^{n\times m}$ is the transpose of $A$.
Furthermore, $I_n\in\R^{n\times n}$ is the identity matrix.

We equip $\R^n$ with the Euclidean norm $\norm{\cdot}$.
Given $\bar x\in\R^n$ and $\delta>0$, $\mathbb B_\delta(\bar x):=\{x\,|\,\norm{x-\bar x}\leq\delta\}$
represents the closed $\delta$-ball around $\bar x$.
For a set $C\subset\R^n$,
$\dist(\bar x,C):=\inf_{x\in C}\norm{x-\bar x}$ is the distance
of $\bar x$ to $C$ with the convention $\dist(\bar x,\emptyset)=\infty$.
For brevity of notation, we make use of
$\bar x+C:=C+\bar x:=\{x+\bar x\,|\,x\in C\}$.

For a lower semicontinuous function $h\colon\R^n\to\overline{\R}$, 
$\dom h:=\{x\in\R^n\,|\,h(x)<\infty\}$ and $\epi h:=\{(x,\alpha)\,|\,h(x)\leq\alpha\}$
are called the domain and epigraph of $h$, respectively,
and we note that $\epi h$ is closed.
Furthermore, $h$ is referred to as proper 
if $\dom h\neq\emptyset$ and $h(x)>-\infty$ is valid for all $x\in\R^n$.
For a nonempty, closed set $C\subset\R^n$, we are concerned with the indicator function
$\delta_C\colon\R^n\to\overline{\R}$ which takes value $0$ on $C$ and value $\infty$
on $\R^n\setminus C$. Obviously, we have $\dom\delta_C=C$ and $\epi\delta_C=C\times[0,\infty)$,
which particularly means that $\delta_C$ is proper and lower semicontinuous.

For a continuously differentiable mapping $\Phi\colon\R^n\to\R^m$ as well as
$\bar x\in\R^n$, $\Phi'(\bar x)\in\R^{m\times n}$ is the Jacobian of $\Phi$ at $\bar x$.
Additionally, for twice continuously differentiable $\varphi\colon\R^n\to\R$,
$\nabla \varphi(\bar x):=\varphi'(\bar x)^\top$ and 
$\nabla^2\varphi(\bar x):=(\nabla \varphi)'(\bar x)$
are used to denote the gradient and the Hessian of $\varphi$ at $\bar x$, respectively.
Partial derivatives w.r.t.\ (with respect to) certain variables are denoted in the usual way.

\subsection{Tangent and normal cones}\label{sec:tangents_and_normals}

For a set $C\subset\R^n$ 
and some point $\bar x\in C$,
we make use of
\begin{align*}
	T_C(\bar x)
	&:=
	\left\{
		d\,\middle|\,\exists d_k\to d,\,t_k\downarrow 0\colon\,\bar x+t_kd_k\in C\,\forall k\in\N
	\right\},
	\\
	T^{\#}_C(\bar x)
	&:=
	\left\{
		d\,\middle|\,\exists x_k\to \bar x,\,d_k\to d\colon\,x_k\in C,\,d_k\in T_C(x_k)\,\forall k\in\N
	\right\},
	\\
	T_C^\textup{P}(\bar x)
	&:=
	\left\{
		d\,|\,\exists x_k\to \bar x,\,d_k\to d,\,t_k\downarrow 0\colon\,x_k\in C,\,x_k+t_kd_k\in C\,\forall k\in\N
	\right\}
\end{align*}
which are referred to as \emph{tangent cone}, \emph{limiting tangent cone}, and \emph{paratingent cone}
to $C$ at $\bar x$, respectively. 
While the tangent cone is a well-known variational object,
the paratingent cone, which seemingly dates back to \cite{Shi1988}, is less popular.
To the best of our knowledge, the limiting tangent cone has been introduced
just recently in \cite[Definition~3.8]{GfrererOutrata2022}.
All these sets are closed cones, obeying the general inclusions
$T_C(\bar x)\subset T_C^{\#}(\bar x)\subset T_C^\textup{P}(\bar x)$.
Note that these inclusions can be strict even for convex polyhedral sets $C$.
Exemplary, we have 
\[
	T_{\R^2_-}(0)=\R^2_-,
	\qquad
	T_{\R^2_-}^{\#}(0)=\{d\,|\,\min(d_1,d_2)\leq 0\},
	\qquad
	T_{\R^2_-}^\textup{P}(0)=\R^2.
\]
The following lemma summarizes some essential calculus rules for the tangent,
limiting tangent, and paratingent cone. 

\begin{lemma}\label{lem:simple_calculus_for_tangents}
	\begin{enumerate}
		\item\label{item:union_rule_for_tangents}
			Let $S_1,\ldots,S_p\subset\R^n$ be closed sets and fix
			$\bar x\in\bigcup_{i=1}^pS_i$. Then
			\[
				\mathcal T_{\bigcup_{i=1}^pS_i}(\bar x)
				\supset
				\bigcup\limits_{i\in I(\bar x)} \mathcal T_{S_i}(\bar x)				
			\]
			holds for each of the tangent cone operators $\mathcal T\in\{T,T^{\#},T^\textup{P}\}$
			where $I(\bar x):=\{i\in\{1,\ldots,p\}\,|\,\bar x\in S_i\}$.	
			Furthermore, for $\mathcal T\in\{T,T^{\#}\}$, we even have equality.		
		\item Let $\Phi\colon\R^n\to\R^m$ be continuously differentiable
			and let $D\subset\R^m$ be a closed set.
			We consider the set $C\subset\R^n$ given by
			\begin{equation}\label{eq:pre-image}
				C:=\{x\,|\,\Phi(x)\in D\}
			\end{equation}
			and fix $\bar x\in C$.
			Then, for each tangent cone operator $\mathcal T\in\{T,T^{\#},T^\textup{P}\}$,
			we always have
			\[
				\mathcal T_C(\bar x)
				\subset
				\{d\,|\,\Phi'(\bar x)\in \mathcal T_D(\Phi(\bar x))\},
			\]
			and equality holds whenever $\Phi'(\bar x)$ possesses full row rank $m$.
		\item Let $\Phi\colon\R^n\to\R^m$ be continuously differentiable
			and let $D\subset\R^m$ be the union of finitely many subspaces of $\R^m$.
			We consider the set $C\subset\R^n$ given in \eqref{eq:pre-image}
			and fix $\bar x\in C$ such that $\Phi'(\bar x)$ possesses full row rank $m$.
			Then we have $T_C(\bar x)=T^{\#}_C(\bar x)$.
			Whenever $D$ is a single subspace, we even have
			$T_C(\bar x)=T^{\#}_C(\bar x)=T^{\textup{P}}_C(\bar x)$.
	\end{enumerate}
\end{lemma}
\begin{proof}
	Let us start with the proof of the first assertion.
	The statement for tangents can be found in \cite[Table~4.1]{AubinFrankowska2009}.
	The proof of the inclusions $\supset$ for the limiting tangent and paratingent
	cone are obvious from $S_j\subset\bigcup_{i=1}^pS_i$ for each $j\in I(\bar x)$.
	For the proof of the converse inclusion for the limiting tangent cone,
	we fix $d\in T^{\#}_{\bigcup_{i=1}^pS_i}(\bar x)$.
	Then we find sequences $x_k\to\bar x$ and $d_k\to d$ such that
	$x_k\in \bigcup_{i=1}^p S_i$ and $d_k\in T_{\bigcup_{i=1}^pS_i}(x_k)=\bigcup_{i\in I(x_k)}T_{S_i}(x_k)$ 
	for each $k\in\N$. Exploiting the pigeonhole principle, we find a set $K\subset\N$ of
	infinite cardinality and some index $j\in\{1,\ldots,p\}$ such that
	$d_k\in T_{S_j}(x_k)$ for each $k\in K$. Closedness of $S_j$ gives
	$j\in I(\bar x)$, and the definition of the limiting tangent cone yields
	$d\in T^{\#}_{S_j}(\bar x)$.
	
	Let us proceed with the proof of the second assertion.
	For tangents, the statement can be found in
	\cite[Table~4.1]{AubinFrankowska2009} and \cite[Exercise~6.7]{RockafellarWets1998}.
	For limiting tangents, this has been proven in \cite[Proposition~4.4]{GfrererOutrata2023}.
	It remains to show the statement for the paratingent cone,
	which can be derived from \cite[Theorem 2.1]{BenkoRockafellar2023}, but we provide a direct proof.\\
	\noindent
	For $d\in T^\textup{P}_C(\bar x)$, 
	there are sequences $x_k\to\bar x$, $t_k\downarrow 0$, and $d_k\to d$
	such that $x_k\in C$ and $x_k+t_kd_k\in C$, 
	i.e., $\Phi(x_k+t_kd_k) = \Phi(x_k)+t_k q_k \in D$ for
	$q_k:=\big(\Phi(x_k+t_kd_k) - \Phi(x_k)\big)/t_k$, for all $k\in\N$.
	Let $L>0$ be some Lipschitz constant of $\Phi'$ around $\bar x$.
	The fundamental theorem of calculus yields
	\[
		\norm{q_k - \Phi'(x_k)d_k}
		=
		\norm{%
		\int_0^1
		\left(\Phi'(x_k + \tau t_k d_k) - \Phi'(x_k)\right)d_k\, \d \tau
		}%
		\leq
		\frac L2 t_k \norm{d_k}^2
	\]
	for each $k\in\N$, and $\Phi'(\bar x)d\in T^\textup{P}_D(\Phi(\bar x))$ follows.\\
	\noindent	
	To prove the converse inclusion, let us assume that $\Phi'(\bar x)$ possesses
	full row rank and pick some $d\in\R^n$ satisfying 
	$\Phi'(\bar x)d\in T^\textup{P}_D(\Phi(\bar x))$.
	Then we find sequences $y_k\to\Phi(\bar x)$, $t_k\downarrow 0$, and $q_k\to \Phi'(\bar x)d$
	such that $y_k\in D$ and $y_k+t_kq_k\in D$ for all $k \in \N$.
	The full rank assumption means that $\Phi$ is so-called \emph{metrically regular} 
	at $(\bar x, \Phi(\bar x))$ with some locally valid regularity constant $\kappa > 0$,
	see \cite[Corollary~3.8]{Mordukhovich2018}.
	Thus, there exists a sequence $(x_k)\subset \R^n$ with $\Phi(x_k) = y_k$ 
	and $\norm{x_k - \bar x} \leq \kappa \norm{y_k - \Phi(\bar x)}$
	for all $k \in \N$, showing $x_k \to \bar x$.
	Using metric regularity of $\Phi$ again yields another sequence $(\hat x_k)\subset \R^n$ 
	with $\Phi(\hat x_k) = y_k+t_kq_k$ (i.e., $\hat x_k \in C$)
	and
	\begin{align*}
		\norm{\hat x_k - (x_k + t_k d)}
		&\leq
		\kappa \norm{y_k+t_kq_k - \Phi(x_k + t_k d)}
		\\
		&=
		t_k \kappa
		\norm{q_k-\frac{\Phi(x_k+t_k d)-\Phi(x_k)}{t_k}}
		\\
		&\leq
		t_k \kappa
		\left(
			\norm{q_k - \Phi'(x_k)d} + \frac L2 t_k \norm{d}^2
		\right)
	\end{align*}
	for all $k \in \N$, where the last inequality follows from the
	same arguments as used to show the first inclusion.
	Consequently, we have $(\hat x_k - x_k)/t_k \to d$, 
	and $d\in T^\textup{P}_C(\bar x)$ follows from
	$x_k + t_k (\hat x_k - x_k)/t_k = \hat x_k \in C$.
	
	The final assertion is a consequence of the first two 
	as the tangent, limiting tangent, and paratingent cone to a subspace 
	$L\subset\R^m$ reduce to the subspace $L$ and, thus, coincide.
\end{proof}

Let us point out that equality does, in general, 
not hold in statement~\ref{item:union_rule_for_tangents}
when considering the paratingent cone to unions of closed sets
even if all involved sets are convex polyhedral.
Exemplary, set $S_1:=\R\times\{0\}$ as well as $S_2:=\{0\}\times\R$ and consider
$\bar x:=(0,0)\in S_1\cup S_2$.
Then $T_{S_1}^\textup{P}(\bar x)\cup T_{S_2}^\textup{P}(\bar x)=S_1\cup S_2$
but $T_{S_1\cup S_2}^\textup{P}(\bar x)=\R^2$.

For a set $C\subset\R^n$
	and some point $\bar x\in C$,
we will exploit
\begin{align*}
	\widehat N_C(\bar x)
	&:=
	\{
		\eta\,|\,\eta^\top (x-\bar x)\leq\oo(\norm{x-\bar x})\,\forall x\in C
	\},
	\\
	N_C(\bar x)
	&:=
	\{
		\eta\,|\,\exists x_k\to\bar x,\,\eta_k\to\eta\colon\,x_k\in C,\,\eta_k\in\widehat N_C(x_k)\,\forall k\in\N
	\},
\end{align*}
the so-called \emph{regular} and \emph{limiting normal cone} to $C$ at $\bar x$.
These sets are closed cones which satisfy $\widehat N_C(\bar x)\subset N_C(\bar x)$
by definition, and $\widehat N_C(\bar x)$ is, additionally, convex.
Whenever $C$ is a convex set, $\widehat N_C(\bar x)$ and $N_C(\bar x)$ are the same
and coincide with the standard normal cone of convex analysis.

\subsection{Generalized derivatives of lower semicontinuous functions}
\label{sec:generalized_differentiation}

Let us consider a proper, lower semicontinuous function $h\colon\R^n\to\overline{\R}$
as well as some point $\bar x\in\dom h$. Then
\begin{align*}
	\partial h(\bar x)
	&:=
	\{z\,|\,(z,-1)\in N_{\epi h}(\bar x,h(\bar x))\},
	\\
	\partial^\infty h(\bar x)
	&:=
	\{z\,|\,(z,0)\in N_{\epi h}(\bar x,h(\bar x))\}
\end{align*}
are called the \emph{limiting} and \emph{singular subdifferential} of $h$ at $\bar x$, respectively. 
It should be noted that the singular subdifferential possesses the equivalent representation
\[
	\partial^\infty h(\bar x)
	=
	\{z\,|\,\exists x_k\to\bar x,\,t_k\downarrow 0,\,z_k\to z\colon\,
		h(x_k)\to h(\bar x),\,z_k\in t_k\partial h(x_k)\,\forall k\in\N\},
\]
while $\partial h(\bar x)$ enjoys the following robustness property:
\[
	\partial h(\bar x)
	=
	\{z\,|\,\exists x_k\to\bar x,\,z_k\to z\colon\,
		h(x_k)\to h(\bar x),\,z_k\in \partial h(x_k)\,\forall k\in\N\}.
\]
If $h$ is convex, $\partial h(\bar x)$ coincides with the standard subdifferential of convex analysis.

For $z\in\R^n$, the function $\mathrm d^2h(\bar x,z)\colon\R^n\to\overline{\R}$ given via
\[
	\mathrm d^2 h(\bar x,z)(d)
	:=
	\liminf\limits_{t\downarrow 0,\,d'\to d}
	\frac{h(\bar x+t d')-h(\bar x)-t\,z^\top d'}{\frac12 t^2} 
\]
is referred to as the \emph{second subderivative} of $h$ at $\bar x$ for $z$.
For more information about this variational tool, 
we refer the interested reader to \cite[Section~13B]{RockafellarWets1998}
and the recent paper \cite{BenkoMehlitz2023} where an overview of available
calculus rules for the second subderivative is provided.

\subsection{Set-valued mappings}\label{sec:set_valued_maps}

\paragraph{Basics}
Throughout the subsection, we consider a set-valued mapping $\Upsilon\colon\R^n\tto\R^m$.
Recall that $\dom\Upsilon:=\{x\,|\,\Upsilon(x)\neq\emptyset\}$, $\gph\Upsilon:=\{(x,y)\,|\,y\in\Upsilon(x)\}$,
and $\ker\Upsilon:=\{x\,|\,0\in\Upsilon(x)\}$
are referred to as the domain, the graph, and the kernel of $\Upsilon$, respectively.
The inverse $\Upsilon^{-1}\colon\R^m\tto\R^n$ of $\Upsilon$ 
is defined via $\gph\Upsilon^{-1}:=\{(y,x)\,|\,y\in\Upsilon(x)\}$. 
Recall that $\Upsilon$ is called locally bounded at $\bar x\in\dom\Upsilon$ 
whenever there exists a neighborhood $U\subset\R^n$ of $\bar x$ such that
$\Upsilon(U):=\bigcup_{x\in U}\Upsilon(x)$ is bounded.

Next, we introduce primal and dual derivatives
of set-valued mappings via tangent and normal cones to their graphs, respectively.
Fix some point $(\bar x,\bar y)\in\gph\Upsilon$.
For each of the tangent cone operators $\mathcal T\in\{T,T^{\#},T^\textup{P}\}$,
the associated derivative $\mathcal D \Upsilon(\bar x,\bar y)\colon\R^n\tto\R^m$,
where $\mathcal D\in\{D,D^{\#},D^\textup{P}\}$ is chosen according to $\mathcal T$,
is defined by
\[
	\gph\mathcal D\Upsilon(\bar x,\bar y)
	:=
	\mathcal T_{\gph\Upsilon}(\bar x,\bar y).
\]
We refer to $D\Upsilon(\bar x,\bar y)$, $D^{\#}\Upsilon(\bar x,\bar y)$, and $D^{\textup{P}}\Upsilon(\bar x,\bar y)$
as the \emph{graphical}, \emph{limiting graphical}, and \emph{strict graphical derivative} 
of $\Upsilon$ at $(\bar x,\bar y)$, respectively. 
For each of the normal cone operators $\mathcal N\in\{\widehat N,N\}$,
the associated derivative $\mathcal D^* \Upsilon(\bar x,\bar y)\colon\R^m\tto\R^n$,
where $\mathcal D^*\in\{\widehat D^*,D^*\}$ is chosen according to $\mathcal N$,
is defined by
\[
	\gph\mathcal D^*\Upsilon(\bar x,\bar y)
	:=
	\{(y^*,x^*)\,|\,(x^*,-y^*)\in \mathcal N_{\gph\Upsilon}(\bar x,\bar y)\}.
\]
We refer to $\widehat D^*\Upsilon(\bar x,\bar y)$ and $D^*\Upsilon(\bar x,\bar y)$
as the \emph{regular} and \emph{limiting coderivative} of $\Upsilon$ at $(\bar x,\bar y)$,
respectively.
In the literature, the strict graphical derivative is often denoted by $D_*\Upsilon(\bar x,\bar y)$.
However, for consistency of our notation and in order to avoid any confusion with the limiting
coderivative, we stick to the notation $D^\textup{P}\Upsilon(\bar x,\bar y)$.

\paragraph{Regularity and Lipschitzian properties}
In this paper, we are concerned with various regularity and
Lipschitzian properties of set-valued mappings.
Let us start by stating the definition of certain regularity properties.
Recall that $\Upsilon$ is said to be \emph{metrically subregular} at $(\bar x,\bar y)\in\gph\Upsilon$
(with constant $\kappa>0$) 
whenever there is a neighborhood $U\subset\R^n$ of $\bar x$ such that
\[
	\forall x\in U\colon\quad
	\dist(x,\Upsilon^{-1}(\bar y))
	\leq
	\dist(\bar y,\Upsilon(x)).
\]
If, additionally, $\bar x$ is an isolated point of $\Upsilon^{-1}(\bar y)$, 
then $\Upsilon$ is said to be \emph{strongly metrically subregular} at $(\bar x,\bar y)$
(with constant $\kappa>0$).
The infimum over all constants $\kappa$ such that the above estimate holds 
is referred to as the modulus of (strong) metric subregularity.
If there is a neighborhood $W\subset\R^n\times\R^m$ of $(\bar x,\bar y)$ such that
$\Upsilon$ is strongly metrically subregular at all points from
$W\cap\gph\Upsilon$ with constant $\kappa>0$, then, according to \cite[Definition~2.7]{GfrererOutrata2023}, 
$\Upsilon$ is called \emph{strongly metrically subregular around} $(\bar x,\bar y)$ 
(with constant $\kappa$), or strongly metrically subregular on a neighborhood of $(\bar x,\bar y)$.
Similarly as above, the infimum over all constants $\kappa$ such that $\Upsilon$ is strongly
metrically subregular around $(\bar x,\bar y)$ with constant $\kappa$ is called the modulus of
strong metric subregularity around $(\bar x,\bar y)$.

Let us proceed with the recapitulation of certain Lipschitzian properties.
We refer to $\Upsilon$ as \emph{calm} at $(\bar x,\bar y)$ (with constant $\kappa>0$)
whenever there is a neighborhood $V\subset\R^m$ of $\bar y$ such that
\[
	\forall x\in\R^n,\,\forall y\in\Upsilon(x)\cap V\colon\quad
	\dist(y,\Upsilon(\bar x)) \leq \kappa\norm{x-\bar x}.
\]
holds. Whenever this estimate can be strengthened to
\[
	\forall x\in\R^n,\,\forall y\in\Upsilon(x)\cap V\colon\quad
	\norm{y-\bar y} \leq \kappa\norm{x-\bar x},
\]
then $\Upsilon$ is said to be \emph{isolatedly calm} at $(\bar x,\bar y)$ (with constant $\kappa>0$),
and we note that this guarantees that $\bar y$ is an isolated point of $\Upsilon(\bar x)$.
Furthermore, $\Upsilon$ is called \emph{isolatedly calm around} $(\bar x,\bar y)$, 
or on a neighborhood of $(\bar x,\bar y)$, (with constant $\kappa>0$) 
if there exists a neighborhood $W\subset\R^n\times\R^m$ of $(\bar x,\bar y)$
such that $\Upsilon$ is isolatedly calm at all points from $W\cap\gph\Upsilon$ with
constant $\kappa$, see \cite[Definition~2.7]{GfrererOutrata2023}.
Similar as above, one can define the moduli of calmness, isolated calmness,
and isolated calmness on a neighborhood.

It is well known that $\Upsilon$ is calm at $(\bar x,\bar y)$
if and only if $\Upsilon^{-1}$ is metrically subregular at $(\bar y,\bar x)$,
see e.g.\ \cite[Theorems~3H.3]{DontchevRockafellar2014}.
Similarly, it is shown in \cite[Theorems~3I.3]{DontchevRockafellar2014}
that $\Upsilon$ is isolatedly calm at $(\bar x,\bar y)$ if and only if
$\Upsilon^{-1}$ is strongly metrically subregular at $(\bar y,\bar x)$.
	The latter always implies and, in the case where $\gph\Upsilon$ is closed locally around $(\bar x,\bar y)$, 
	is even equivalent to the so-called \emph{Levy--Rockafellar criterion}
	\[
		\ker D\Upsilon^{-1}(\bar y,\bar x) = \{0\},
	\]
	see \cite[Theorem~4E.1]{DontchevRockafellar2014},
	and this, in turn, is obviously the same as
	\begin{equation}\label{eq:LevyRockafellar_criterion}
		D\Upsilon(\bar x,\bar y)(0)=\{0\}.
	\end{equation}
Furthermore, it should be noted that $\Upsilon$ is isolatedly calm around $(\bar x,\bar y)$
if and only if $\Upsilon^{-1}$ is strongly metrically subregular around $(\bar y,\bar x)$,
see \cite[Lemma~2.8]{GfrererOutrata2023}.
By \cite[Theorem~2.9(iii)]{GfrererOutrata2023}, whenever $\gph\Upsilon$ is closed
locally around $(\bar x,\bar y)$, the latter is equivalent to
\[
	\ker D^{\#}\Upsilon^{-1}(\bar y,\bar x) = \{0\},
\]
and this is the same as
\begin{equation}\label{eq:GfrererOutrata_criterion}
	D^{\#}\Upsilon(\bar x,\bar y)(0) = \{0\}.
\end{equation}
	Again, it is easy to see that the necessity of this criterion does not require any closedness of $\gph\Upsilon$.
We would like to recall that $\Upsilon$ possesses the \emph{Aubin property} at $(\bar x,\bar y)$,
see \cite[Section~3E]{DontchevRockafellar2014} for the definition of this well-established property,
if and only if the so-called \emph{Mordukhovich criterion}
\begin{equation}\label{eq:Mordukhovich_criterion}
	D^*\Upsilon(\bar x,\bar y)(0) = \{0\}
\end{equation}
is valid
	and $\gph\Upsilon$ is closed locally around $(\bar x,\bar y)$,
see e.g.\ \cite[Theorem~3.3]{Mordukhovich2018}.
Finally, following \cite[Theorems~2.6, 2.7(iii)]{GfrererOutrata2022},
$\Upsilon$ possesses a \emph{single-valued Lipschitz continuous localization} around $(\bar x,\bar y)$,
i.e., there are neighborhoods $U\subset\R^n$ of $\bar x$ and $V\subset\R^m$ of $\bar y$ as well as a 
locally Lipschitz continuous function $\upsilon\colon U\to\R^m$ with $\upsilon(\bar x)=\bar y$
and $(U\times V)\cap\gph\Upsilon=\{(x,\upsilon(x))\,|\,x\in U\}$, if and only if
\begin{equation}\label{eq:Strict_graph_der_criterion}
	D^\textup{P}\Upsilon(\bar x,\bar y)(0)=\{0\}
\end{equation}
and \eqref{eq:Mordukhovich_criterion} hold simultaneously.
	Clearly, whenever $\Upsilon$ possesses a single-valued Lipschitz continuous localization around $(\bar x,\bar y)$,
	then $\gph\Upsilon$ is trivially closed locally around $(\bar x,\bar y)$.
	Conversely, \eqref{eq:Strict_graph_der_criterion} also implies this local closedness of $\gph\Upsilon$
	around $(\bar x,\bar y)$ which can be distilled from \cite[formula (2.7)]{BenkoRockafellar2023}. 

\paragraph{Inner calmness*}

While the aforementioned regularity and Lipschitzian notions are comparatively common,
we will now recall the concept of so-called \emph{inner calmness* in the fuzzy sense} which
has been recently introduced in \cite[Definition~2.6]{Benko2021}.

\begin{definition}\label{def:IC*_in_fuzzy_sense}
	Let $\Upsilon\colon\R^n\tto\R^m$ be a set-valued mapping,
	and fix $\bar x\in\dom\Upsilon$.
	For some direction $d\in\R^n$,
	we say that $\Upsilon$ is \emph{inner calm* in the fuzzy sense}
	w.r.t.\ $\dom\Upsilon$ at $\bar x$ in direction $d$ 
	if either $d\notin T_{\dom\Upsilon}(\bar x)$
	or if there is a constant $\kappa_d>0$ such that we find $\bar y\in\R^m$ as well as
	sequences $d_k\to d$, $t_k\downarrow 0$, and $y_k\to\bar y$
	such that $\bar x+t_kd_k\in\dom\Upsilon$, $y_k\in\Upsilon(\bar x+t_kd_k)$,
	and $\norm{y_k-\bar y}\leq\kappa_dt_k\norm{d_k}$ for all $k\in\N$.
	The infimum over all constants $\kappa_d$ with this property is called
	the modulus of inner calmness* in the fuzzy sense w.r.t.\ $\dom\Upsilon$
	at $\bar x$ in direction $d$, and it is set to $0$ 
	if $d\notin T_{\dom\Upsilon}(\bar x)$.
	Furthermore, $\Upsilon$ is called inner calm* in the fuzzy sense 
	w.r.t.\ $\dom\Upsilon$ at $\bar x$ whenever it is inner calm* in the fuzzy
	sense w.r.t.\ $\dom\Upsilon$ at $\bar x$ in each unit direction.
\end{definition}

Inner calmness* in the fuzzy sense was introduced for its role in certain calculus rules,
see \cite[Theorem 4.1]{Benko2021} and \cite[Theorem 3.1]{BenkoMehlitz2022}.
In \cite[Theorem 4.1]{BenkoMehlitz2022}, 
it was shown that this notion is also necessary for validity of those calculus rules.
Perhaps more importantly, a (Lagrange) multiplier mapping 
associated with NLPs always has this property, 
see \cite[Theorem 3.9]{Benko2021}.
Both these basic results will be essential for this paper.

Let us mention the simpler concept of inner calmness*, 
see \cite[Definition~2.2]{Benko2021}.
This notion arises naturally from the known concept of inner semicompactness, 
see e.g.\ \cite[Section~2]{Benko2021} for a definition, as it simply enriches it with a rate.
By \cite[Theorem~3.4]{Benko2021}, polyhedral mappings, i.e., mappings whose graphs
can be represented as the union of finitely many convex polyhedral sets,
enjoy the inner calmness* property w.r.t.\ their domain at each point of their domain.
Moreover, this notion proved to be essential 
for the calculus of second subderivatives in \cite{BenkoMehlitz2023}.
The problem is, however, that the aforementioned multiplier mapping 
does not seem to satisfy inner calmness* automatically
- that is why the weaker fuzzy version from \cref{def:IC*_in_fuzzy_sense} has been introduced.

Below, we introduce a locally uniform version of inner calmness* in the fuzzy sense.

\begin{definition}\label{def:IC*_in_fuzzy_sense_related_notions}
	Let $\Upsilon\colon\R^n\tto\R^m$ be a set-valued mapping,
	and fix $\bar x\in\dom\Upsilon$.
	For some direction $d\in\R^n$,
	we refer to $\Upsilon$ as \emph{locally uniformly inner calm* in the fuzzy sense}
	w.r.t.\ $\dom\Upsilon$ at $\bar x$ in direction $d\in\R^n$ 
	if there are neighborhoods $U\subset\R^n$ of
	$\bar x$ and $N \subset \R^n$ of $d$ as well as a constant $\kappa>0$ such that
	$\Upsilon$ is inner calm* in the fuzzy sense w.r.t.\ $\dom\Upsilon$ at each
	point $x\in U\cap \dom\Upsilon$ in each direction $d' \in N$ 
	with a modulus not larger than $\kappa$.
	The infimum over all constants $\kappa_d$ with this property is called
	the modulus of locally uniform inner calmness* in the fuzzy sense
	w.r.t.\ $\dom\Upsilon$ at $\bar x$ in direction $d$.
	Furthermore, $\Upsilon$ is called locally uniformly inner calm* in the fuzzy sense
	w.r.t.\ $\dom\Upsilon$ at $\bar x$ 
	whenever it is locally uniformly inner calm* in the fuzzy sense
	w.r.t.\ $\dom\Upsilon$ at $\bar x$ in each unit direction.
\end{definition}

\begin{remark}\label{Rem:uniform_ic*_automatic}
If, in the setting of \cref{def:IC*_in_fuzzy_sense_related_notions}, 
$d \notin T^{\#}_{\dom\Upsilon}(\bar x)$,
then there exist neighborhoods $U\subset\R^n$ of $\bar x$ and $N \subset \R^n$ of $d$
such that 
\[
	\{(x',d')\,|\,x'\in\dom\Upsilon,\,d'\in T_{\dom\Upsilon}(x')\}\cap (U\times N) 
	= 
	\emptyset, 
\]
and so, for every $x'\in U\cap \dom\Upsilon$, we have $T_{\dom\Upsilon}(x') \cap N = \emptyset$.
Consequently, $\Upsilon$ is locally uniformly inner calm* in the fuzzy sense 
w.r.t.\ $\dom\Upsilon$ at $\bar x$ in direction $d$ with modulus $0$.
\end{remark}

It should be mentioned that \cite[Theorem~3.4]{Benko2021} actually yields that
polyhedral mappings enjoy also the locally uniform fuzzy inner calmness* property 
w.r.t.\ their domain at each point of their domain.
In \cite[Section~4]{BenkoMehlitz2022}, 
the interested reader can find an overview of various sufficient conditions for
inner calmness* in the fuzzy sense and related notions.
Exemplary, the following lemma provides sufficient conditions for the presence of inner calmness* 
in the fuzzy sense and its locally uniform version in terms of the primal
derivative criteria \eqref{eq:LevyRockafellar_criterion} and \eqref{eq:GfrererOutrata_criterion}.

\begin{lemma}\label{lem:sufficient_conditions_IC*}
Let $\Upsilon\colon\R^n\tto\R^m$ be a set-valued mapping with a closed graph.
Fix $(\bar x,\bar y)\in\gph\Upsilon$.
Let $\Upsilon$ be locally bounded at $\bar x$
and assume that $\Upsilon$ possesses convex images locally around $\bar x$.
If \eqref{eq:LevyRockafellar_criterion} (or even \eqref{eq:GfrererOutrata_criterion}) holds,
then $\Upsilon$ is (locally uniformly) inner calm* in the fuzzy sense w.r.t.\ $\dom\Upsilon$ at $\bar x$.	
\end{lemma}
\begin{proof}
	Note that \eqref{eq:LevyRockafellar_criterion} implies $\Upsilon(\bar x) = \{\bar y\}$
	since the local isolatedness of $\bar y$ becomes global 
	due to convexity of the set $\Upsilon(\bar x)$.
	The first assertion, thus, follows from \cite[Lemma~4.3(ii)]{BenkoMehlitz2022}
	since the assumed local boundedness of $\Upsilon$ at $\bar x$ shows that
	$\Upsilon$ is inner semicompact w.r.t.\ its domain at $\bar x$.
	
	For the proof of the second assertion, 
	let us first note that \eqref{eq:GfrererOutrata_criterion} yields
	isolated calmness of $\Upsilon$ on a neighborhood of $(\bar x,\bar y)$,
	i.e., at all points from a neighborhood of $(\bar x,\bar y)$ from the
	graph of $\Upsilon$, this mapping is isolatedly calm with a
	uniform modulus $\kappa>0$.
	Due to \cite[Theorem~4E.1]{DontchevRockafellar2014},
	this shows that
	\[
		d^y\in D\Upsilon(x,y)(d^x)
		\quad\Longrightarrow\quad
		\norm{d^y}\leq\kappa\norm{d^x}
	\]
	is valid for all $(d^x,d^y),(x,y)\in\R^n\times\R^m$ such that
	$(x,y)\in\gph\Upsilon$ is close enough to $(\bar x,\bar y)$.
	Now, the statement follows from \cite[Theorem~4.1]{BenkoMehlitz2022},
	noting that, due to \cite[Lemma~2.1]{Benko2021}, $\dom\Upsilon$ is closed 
	locally around $\bar x$ since $\Upsilon$ is locally bounded at this point.
\end{proof}

Let us briefly mention that the assumptions of \cref{lem:sufficient_conditions_IC*}
imply much stronger properties.
For instance, assuming \eqref{eq:LevyRockafellar_criterion},
not only we can drop ``in the fuzzy sense'' 
(as clear already from \cite[Lemma~4.3(ii)]{BenkoMehlitz2022}),
but due to $\Upsilon(\bar x) = \{\bar y\}$, 
we can actually apply \cite[Lemma~4.3(i)]{BenkoMehlitz2022}
to also drop the star, arriving simply at inner calmness of $\Upsilon$, 
see \cite{BenkoMehlitz2022} for the definition and more details.
For the purposes of this paper, however, \cref{lem:sufficient_conditions_IC*} will suffice.

The final lemma of this subsection, 
which has been motivated by \cite[Theorem~3.1]{BenkoMehlitz2022},
addresses a calculus rule for tangent cones to the domain of a given set-valued mapping.
Here, inner calmness* in the fuzzy sense is of special importance.

\begin{lemma}\label{lem:calculus_for_limiting_tangent_cone}
	Let $\Upsilon\colon\R^n\tto\R^m$ be a set-valued mapping with a closed graph
	and fix $\bar x\in\dom\Upsilon$.
	For $\mathcal T := T$ and $\mathcal D := D$ 
	($\mathcal T := T^{\#}$ and $\mathcal D := D^{\#}$),
	the following assertions hold.
	\begin{enumerate}
		\item\label{a} We have
			\[
				d^x \in \mathcal T_{\dom\Upsilon}(\bar x)
				\quad \Longleftarrow \quad
				d^x \in \bigcup\limits_{\bar y\in\Upsilon(\bar x)}\dom \mathcal D\Upsilon(\bar x,\bar y),
			\]
			and the converse implication holds whenever $\Upsilon$ is
				locally bounded at $\bar x$ and
			(locally uniformly) inner calm* in the fuzzy sense w.r.t.\ $\dom\Upsilon$ at $\bar x$ in direction $d^x$.
		\item\label{b} We have
			\[
				\mathcal T_{\dom\Upsilon}(\bar x)
				\supset
				\bigcup\limits_{\bar y\in\Upsilon(\bar x)}\dom \mathcal D\Upsilon(\bar x,\bar y),
			\]
			and the converse inclusion holds whenever $\Upsilon$ is
				locally bounded at $\bar x$ and
			(locally uniformly) inner calm* in the fuzzy sense w.r.t.\ $\dom\Upsilon$ at $\bar x$.
	\end{enumerate}	
\end{lemma}
\begin{proof}
	Let us start to prove statement~\ref{a} for $\mathcal T := T$ and $\mathcal D := D$.	
	The general implication trivially follows by definition of the graphical derivative.
	The converse implication is a consequence of \cite[Theorems~3.1 and~4.1]{BenkoMehlitz2022} 
	since $\dom\Upsilon$ is locally closed around $\bar x$.
	The latter can be distilled from \cite[Lemma~2.1]{Benko2021}
	as $\Upsilon$ is assumed to be locally bounded at $\bar x$.

	For the proof of statement~\ref{a} for $\mathcal T := T^{\#}$ and $\mathcal D := D^{\#}$,
	fix $\bar y\in\Upsilon(\bar x)$ and $d^x\in\dom D^{\#}\Upsilon(\bar x,\bar y)$.
	Then we find $d^y\in\R^m$ such that $(d^x,d^y)\in T^{\#}_{\gph\Upsilon}(\bar x,\bar y)$.
	By definition of the limiting tangent cone,
	there are sequences $x_k\to\bar x$, $y_k\to\bar y$, $d^x_k\to d^x$, 
	and $d^y_k\to d^y$ such that
	$(x_k,y_k)\in\gph \Upsilon$ and $(d^x_k,d^y_k)\in T_{\gph\Upsilon}(x_k,y_k)$ 
	for each $k\in\N$.
	The first assertion gives $d^x_k\in T_{\dom\Upsilon}(x_k)$ for large enough $k\in\N$, 
	and the definition of the limiting tangent cone yields 
	$d^x\in T^{\#}_{\dom\Upsilon}(\bar x)$.

	To show the converse implication, let us assume that $\Upsilon$ is
	locally bounded at $\bar x$ and
	locally uniformly inner calm* in the fuzzy sense 
	w.r.t.\ $\dom\Upsilon$ at $\bar x$ in direction $d^x\in T^{\#}_{\dom\Upsilon}(\bar x)$ with modulus $\kappa>0$.
	Consider sequences $x_k\to\bar x$ and $d^x_k\to d^x$
	with $x_k\in\dom\Upsilon$ and $d^x_k \in T_{\dom\Upsilon}(x_k)$ for each $k \in \N$.
	We now apply \cite[Theorem~3.1]{BenkoMehlitz2022} 
	in order to find sequences $(y_k)$ and $(d^y_k)$ satisfying $y_k \in \Upsilon(x_k)$,
	$d^y_k \in D\Upsilon(x_k,y_k)(d^x_k)$, and $\nnorm{d^y_k} \leq \kappa \nnorm{d^x_k}$ 
	for all sufficiently large $k\in\N$.
	The local boundedness of $\Upsilon$ at $\bar x$ yields, 
	by passing to a subsequence if necessary,
	$y_k\to\bar y$ for some $\bar y\in\R^m$,
	and as $(d^y_k)$ is bounded due to boundedness of $(d^x_k)$,
	there is some $d^y\in\R^m$ such that $d^y_k\to d^y$
	holds along a subsequence if necessary.
	Finally, $\bar y \in \Upsilon(\bar x)$ is obtained from the closedness of $\gph\Upsilon$
	and $(d^x,d^y)\in T^{\#}_{\gph\Upsilon}(\bar x,\bar y)$, i.e.,
	$d^y\in D^{\#}\Upsilon(\bar x,\bar y)(d^x)$, 
	follows by the definition of the limiting tangent cone and the associated derivative.

	Statement~\ref{b} clearly follows from~\ref{a}.
\end{proof}

\subsection{On the semismooth* Newton method}\label{sec:ss*}

In this subsection, we briefly recall the semismooth* Newton method from \cite{GfrererOutrata2021}
for solving the generalized equation
\begin{equation}\label{eq:GE}
    0 \in \Upsilon(x),
\end{equation}
where $\Upsilon\colon \R^n \tto \R^n$ is a given set-valued mapping with a closed graph.
We will only detail those parts of the method which are of interest in this paper,
and we refer to \cite{GfrererOutrata2021} 
as well as \cite{GfrererOutrata2022,GfrererOutrataValdman2022} 
for more details about the algorithm and its extensions.

In order to outline the method, we will need some notation.
Given $(x,y) \in \gph \Upsilon$, we denote by $\mathcal{A}\Upsilon(x,y)$ 
the collection of all pairs of
matrices $(A,B)\in\R^{n\times n}\times\R^{n\times n}$ such that there are $n$ pairs
$(y_i^*,x_i^*) \in \gph D^*\Upsilon(x,y)$, $i = 1,\ldots,n$, 
and the $i$-th row of $A$ and $B$ are
$(x_i^*)^\top$ and $(y_i^*)^\top$, respectively. Furthermore, we make use of
\[
	\mathcal{A}_{\textrm{reg}}\Upsilon(x,y) 
	:= 
	\{(A,B) \in \mathcal{A}\Upsilon(x,y) \,|\, A \textrm{ is an invertible matrix}\}.
\]
According to \cite[Algorithm 2]{GfrererOutrata2021}, 
the semismooth* Newton method can be stated as follows.

\begin{algorithm}[Semismooth* Newton method]\leavevmode
	\label{alg:SSN}
	\begin{algorithmic}[1]
		\REQUIRE Starting point $x_0\in\R^n$
		\STATE Set $k := 0$.
		\WHILE{$0\notin\Upsilon(x_k)$}
		\STATE\label{step:approximation_step}
			Compute $(\hat x_k,\hat y_k)\in\gph\Upsilon$ close to $(x_k,0)$ such that
			$\mathcal A_\textup{reg}\Upsilon(\hat x_k,\hat y_k)\neq\emptyset$. 
		\STATE\label{step:Newton_step} 
			Select $(A_k,B_k)\in\mathcal A_\textup{reg}\Upsilon(\hat x_k,\hat y_k)$ and set
			$x_{k+1}:=\hat x_k-A^{-1}_kB_k\hat y_k$.
      	\STATE Set $ k \leftarrow k + 1 $.
		\ENDWHILE
		\RETURN $x_k$
	\end{algorithmic}
\end{algorithm}

Let us briefly mention that \cref{step:approximation_step,step:Newton_step} of \cref{alg:SSN} 
are referred to as \emph{approximation} and \emph{Newton step} in the literature, respectively.

We want to focus on regularity assumptions
that yield superlinear convergence of a sequence generated by \cref{alg:SSN}.
In \cite[Theorem~4.7]{GfrererOutrata2021}, the authors
show that semismoothness* and strong metric regularity,
see \cite[Section~3G]{DontchevRockafellar2014} for a precise definition, 
are sufficient for that purpose.
For brevity of presentation, we are not going to state the definition of semismoothness*
but refer the interested reader to \cite[Section~3]{GfrererOutrata2021} where this property
is introduced and several sufficient criteria for its validity can be found.
Here, let us just mention that set-valued mappings whose graphs are closed subanalytic sets
or can be represented as the union of finitely many closed convex sets 
are semismooth*, see \cite[Proposition~2.10]{GfrererOutrata2022}, 
so this class of mappings is rather large.
Later, in \cite{GfrererOutrata2022}, the same authors
design a so-called SCD version of this method,
where SCD abbreviates \emph{subspace containing derivative},
and show in \cite[Proposition~5.5, Corollary 5.6]{GfrererOutrata2022} that 
so-called SCD semismoothness* and SCD regularity 
are sufficient for its local superlinear convergence.
Again, we are not going to present the definition of SCD mappings and SCD regularity.
A detailed introduction to these concepts can be found in 
\cite[Sections~3 and 4]{GfrererOutrata2022}.
SCD regularity is often equivalent to and always implied
by strong metric subregularity on a neighborhood, see \cite[Theorem~6.2]{GfrererOutrata2022}.
We will now show that strong metric subregularity on a neighborhood is actually sufficient for 
the desired local superlinear convergence properties of \cref{alg:SSN}.
To this end, we first prove that the essential
result \cite[Theorem~4.1]{GfrererOutrata2021}
remains true with strong metric regularity
replaced by strong metric subregularity.

\begin{proposition}\label{prop:SSN_approximation_step}
 	Let $\Upsilon\colon\R^n\tto\R^n$ be a set-valued mapping with a closed graph,
 	and assume that $\Upsilon$ is strongly metrically subregular 
 	at $(\hat x, \hat y)\in\gph\Upsilon$ with modulus $\kappa > 0$.
 	Then there is a matrix $B\in\R^{n\times n}$ with $\norm{B} \leq \kappa$ such that
 	$(I_n,B) \in \mathcal{A}_{\textup{reg}}\Upsilon(\hat x, \hat y)$.
 	Particularly, $\mathcal A_\textup{reg}\Upsilon(\hat x,\hat y)$ is nonempty.
\end{proposition}
\begin{proof}
 	Strong metric subregularity of $\Upsilon$ at $(\hat x,\hat y)$
 	particularly means that $\hat x$ is an isolated point of $\Upsilon^{-1}(\hat y)$, 
 	which implies $N_{\Upsilon^{-1}(\hat y)}(\hat x) = \R^n$.
 	Hence, for each $i=1,\ldots,n$, we have $\mathtt e_i\in N_{\Upsilon^{-1}(\bar y)}(\hat x)$,
 	so that metric subregularity of $\Upsilon$ at $(\hat x,\hat y)$
 	with modulus $\kappa$ and \cite[Theorem~3.2]{BenkoMehlitz2022} yield 
 	the existence of $y^*_i\in\R^n$ with
 	$\norm{y^*_i}\leq\kappa\norm{\mathtt e_i}=\kappa$
 	and
 	$\mathtt e_i\in D^*\Upsilon(\hat x,\hat y)(y^*_i)$.
 	Choosing a suitable matrix norm, the claim follows by definition of
 	$\mathcal A_\textup{reg}\Upsilon(\hat x,\hat y)$.
\end{proof}

With this at hand, the following convergence result
follows by analogous arguments as used to prove
\cite[Theorems~4.4 and 4.7]{GfrererOutrata2021}.

\begin{theorem}\label{The:ss*_superlinear_convergence}
	Let $\Upsilon\colon\R^n\tto\R^n$ be a set-valued mapping with a closed graph.
	Assume that $\Upsilon$ is semismooth* at $(\bar x, 0) \in  \gph \Upsilon$,
	and that there are constants $L, C  > 0$ such that, for every
	$x \notin \Upsilon^{-1}(0)$ sufficiently close to $\bar x$, we have
	$\mathcal{G}^{L,C}_{\Upsilon,\bar x}(x) \neq \emptyset$, where
	\[
		\mathcal{G}^{L,C}_{\Upsilon,\bar x}(x)
		:=
		\left\{
		((\hat x, \hat y), (A,B)) \in \gph \Upsilon \times \mathcal{A}_{\textup{reg}}\Upsilon(\hat{x},\hat{y}) 
		\,\middle|\,
		\begin{aligned}
			&\nnorm{(\hat{x}- \xb,\hat{y})} \leq L \norm{x - \xb},
			\\
			&\nnorm{A^{-1}}\norm{[A \,\vert\, B]}_F \leq C
		\end{aligned}
		\right\}
	\]	
	and $\norm{\cdot}_F$ denotes the Frobenius norm.
	The latter holds true with $L := \beta + 1$ and
	$C := \sqrt{n(1+\kappa^2)}$ for any $\beta \geq 1$
	provided $\Upsilon$ is strongly metrically subregular around $(\xb,0)$
	with modulus $\kappa$.
	Then there exists some $\delta > 0$ such that, for each starting point
	$x_0 \in  \B_\delta(\xb)$, \cref{alg:SSN} either stops
	after finitely many iterations at a solution of \eqref{eq:GE} or
	produces a sequence $(x_k)$ which converges superlinearly to $\xb$,
	provided we choose a quadruple
	$((\hat x_k, \hat y_k), (A_k,B_k)) \in  \mathcal{G}^{L,C}_{\Upsilon,\bar x}(x_k)$
	in each iteration $k\in\N$.
\end{theorem}

\begin{remark}
	Apart from a regularity assumption, we also need the underlying set-valued mapping 
	to be semismooth* in order to show superlinear convergence of \cref{alg:SSN}.
	As mentioned above, semismoothness* is very often present in applications.
	In this paper, we thus pay no attention to this assumption 
	and focus merely on the issue of regularity.
	Similarly, we do not discuss how to carry out \cref{step:approximation_step}
	or how to choose the pair of matrices $(A_k,B_k)$ in \cref{step:Newton_step} 
	and refer the interested reader to
	\cite{GfrererOutrata2021,GfrererOutrata2022,GfrererOutrataValdman2022}.
	Let us just mention that it is actually not necessary 
	to compute the coderivative of the mapping $\Upsilon$,
	which may be quite difficult. 
	The semismooth* Newton method from \cref{alg:SSN} and its extensions
	offer plenty of flexibility regarding the implementation 
	of the approximation and Newton step,
	and the above papers provide a lot of details.
\end{remark}

\section{Perturbing optimization problems}\label{sec:perturbed_optimization_problems}

In this section, we discuss the model problem 
and provide a comprehensive analysis of some associated perturbation mappings.
Particularly, we present the key result \cref{prop:IC_M} which is the foundation of the paper.

\subsection{The model problem and associated perturbation mappings}\label{sec:model}

For easier orientation, we recall here our model mentioned in \cref{sec:introduction}.
The target problem \eqref{eq:Pbar} is given as
\begin{equation*}
	\min\limits_x\quad f(x,0) \quad\text{s.t.}\quad x\in\R^n 
\end{equation*}
for a proper, lower semicontinuous function $f\colon\R^n\times\R^m\to\overline{\R}$.
The latter is embedded in a parametrized family of problems \eqref{eq:P(u,v)} given by
\begin{equation*}
	\min\limits_x\quad f(x,u) - v ^\top x \quad\text{s.t.}\quad x\in\R^n
\end{equation*}
involving parameters $v\in\R^n$ and $u\in\R^m$.
Hence, \eqref{eq:Pbar} and \hyperref[eq:P(u,v)]{(P$(0,0)$)} are the same.
Subsequently, we use $(\bar v,\bar u):=(0,0)$ to denote the pair of reference 
parameters so that \hyperref[eq:P(u,v)]{(P$(\bar v,\bar u)$)}
recovers \eqref{eq:Pbar}.

If $x\in\R^n$ is locally optimal for the perturbed problem \eqref{eq:P(u,v)}, 
a suitable first-order optimality condition simply
reads $v \in \partial_x f(x,u)$, where $\partial_x f(x,u)$ denotes the subdifferential of the function $f(\cdot,u)$ at $x$.
Each point which satisfies this first-order condition will be referred to as \emph{stationary} for \eqref{eq:P(u,v)}.	
Throughout the section, $\bar x\in\R^n$ is a stationary point of \eqref{eq:Pbar}.

	In order to introduce multipliers, we define the mapping
	$Y\colon\R^n\times\R^m\times\R^n\tto\R^m$ associated with \eqref{eq:P(u,v)} given by
	\[
		Y(v,u,x):=\{y\,|\, (v,y) \in \partial f(x,u)\}.
	\]
	Properties of mapping $Y$ will play an important role in the whole paper.
	Desirable behavior of $Y$ is typically secured by the basic qualification condition
	\begin{equation}\label{eq:basic_CQ}
		(0,y)\in\partial^\infty f(\bar x,\bar u)\quad\Longrightarrow\quad y=0,
	\end{equation}
	but we do not make it a standing assumption for reasons explained in the next section.

	We are mainly interested in the case where $f$ is of the composite form \eqref{eq:GNLP_setting}.
	\begin{assumption}\label{ass:composite_setting}
		The function $f$ is given in the form \eqref{eq:GNLP_setting}, i.e.,
		\begin{equation*}
			f(x,u) := f_0(x) + g(F(x) + u)
		\end{equation*}
		for twice continuously differentiable functions $f_0\colon\R^n\to\R$ and $F\colon\R^n\to\R^m$
		as well as a proper, lower semicontinuous, convex function $g\colon\R^m\to\overline{\R}$.
	\end{assumption}
	Whenever \cref{ass:composite_setting} is exploited, we will mention this clearly.
	In the presence of \cref{ass:composite_setting}, the qualification condition \eqref{eq:basic_CQ} reads
		\begin{equation}\label{eq:basic_CQ_composite}
			F'(\bar x)^\top y=0,\,y\in\partial^\infty g(F(\bar x))
			\quad\Longrightarrow\quad
			y=0,
		\end{equation}
	while the (multiplier-dependent) optimality condition $(v,y) \in \partial f(x,u)$ 
	(i.e., $y \in Y(v,u,x)$) takes the form
	\begin{equation}\label{eq:perturbed_stationarity_system}
		\nabla_xL(x,y)=v,\qquad
		y\in\partial g(F(x) + u),
	\end{equation}
	where $L\colon\R^n\times\R^m\to\R$ denotes the \emph{Lagrangian function} given by
	\[
		L(x,y):=f_0(x)+y^\top F(x).
	\]
	Let us now sum up the consequences of such composite structure for the mapping $Y$.
	\begin{proposition}\label{Pro:Basic_properties_from_assumptions}
		Let \cref{ass:composite_setting} hold.
		Then $Y$ possesses the following basic properties.
		\begin{enumerate}
			\item\label{item:Y_convex_valued} 
				For every $(v,u,x) \in \dom Y$, we have $v \in \partial_x f(x,u)$
				and $Y(v,u,x)$ is convex.
			\item\label{item:Y_gph_closed} 
				The set $\gph Y$ is closed.
		\end{enumerate}
		If the qualification condition \eqref{eq:basic_CQ_composite} holds at $\bar x$,
		then there is a closed neighborhood $W\subset\R^n\times\R^m\times\R^n$ of $(\bar v, \bar u, \bar x)$ 
		such that $Y$ satisfies the following assertions as well. 
		\begin{itemize}
		\item[(c)]\label{item:Y_dom_closed} 
		The set $W\cap\dom Y$ is closed and, for all $(v,u,x) \in W$, we have
		\begin{equation}\label{eq:robust_multiplier_rule}
			\begin{aligned}
			v \in \partial_x f(x,u)
			\quad &\Longleftrightarrow \quad
			\exists\, y\in\R^m\colon\, (v,y) \in \partial f(x,u)
			\\
			\quad &\Longleftrightarrow \quad
			(v,u,x) \in \dom Y.
			\end{aligned}
		\end{equation}
		Particularly, $\gph \partial_x f$ is also closed 
		locally around $((\bar x,\bar u),\bar v)$.
		\item[(d)]\label{item:Y_loc_bounded} 
			The set $Y(W) \subset \R^m$ is bounded.
		\end{itemize}
	\end{proposition}
	\begin{proof}
	Statement~\ref{item:Y_convex_valued} is a simple consequence of convexity of $g$ and \cite[Corollary~10.11]{RockafellarWets1998}.
	To show statement~\ref{item:Y_gph_closed}, suppose
	there are sequences $(v_k)$, $(u_k)$, $(x_k)$, and $(y_k)$ with
	$((v_k,u_k,x_k),y_k) \in \gph Y$ for each $k \in \N$
	and a quadruple $(\hat v,\hat u,\hat x, \hat y)\in\R^n\times\R^m\times\R^n\times\R^m$
	such that $(v_k,u_k,x_k,y_k) \to (\hat v,\hat u,\hat x, \hat y)$.
	By \eqref{eq:perturbed_stationarity_system}, we easily get
	$\hat v = \nabla_xL(\hat x,\hat y)$ in the limit by continuous differentiability of $f_0$ and $F$ as well as
	$y_k \in \partial g(F(x_k) + u_k)$ for all $k \in \N$ with $F(x_k) + u_k \to F(\hat x) + \hat u$.
	On the one hand, we get $\liminf_{k \to \infty} g(F(x_k) + u_k) \geq g(F(\hat x) + \hat u)$ by lower semicontinuity of $g$.
	On the other hand, for each $k \in \N$, we have
	\[
		g(F(\hat x) + \hat u) \geq g(F(x_k) + u_k) + y_k ^\top (F(\hat x) + \hat u - F(x_k) - u_k)
	\]
	by convexity of $g$,
	yielding $g(F(\hat x)+\hat u)\geq\limsup_{k\to\infty}g(F(x_k)+u_k)$.
	Thus, the convergence $g(F(x_k) + u_k) \to g(F(\hat x) + \hat u)$ follows and so does $\hat y \in \partial g(F(\hat x) + \hat u)$ in turn.
	This yields $(\hat v,\hat u,\hat x, \hat y) \in \gph Y$ and proves closedness of $\gph Y$, i.e., assertion~\ref{item:Y_gph_closed}.

	Having also the qualification condition \eqref{eq:basic_CQ_composite},
	\cite[Proposition~2.2]{LevyPoliquinRockafellar2000} yields that
	$f$ satisfies the so-called parametric version of \emph{continuous prox-regularity} 
	defined in \cite[Definition~2.1]{LevyPoliquinRockafellar2000}.
	The rest will now follow utilizing results from \cite{LevyPoliquinRockafellar2000}.

	Let us now prove assertion \hyperref[item:Y_dom_closed]{(c)}. 
	The first equivalence in \eqref{eq:robust_multiplier_rule} is stated in
	\cite[Proposition~3.4]{LevyPoliquinRockafellar2000},
	while the second follows by definition of $Y$.
	The local closedness of $\dom Y$ thus follows 
	from the local closedness of $\gph \partial_x f$ stated in
	\cite[Proposition~3.2(b)]{LevyPoliquinRockafellar2000}.

	To show assertion \hyperref[item:Y_loc_bounded]{(d)}, 
	suppose that there are sequences $(v_k)$, $(u_k)$, $(x_k)$, and $(y_k)$ with
	$y_k \in Y(v_k,u_k,x_k)$ for each $k \in \N$, $(v_k,u_k,x_k) \to (\bar v,\bar u,\bar x)$, 
	and $\norm{y_k} \to \infty$.
	Thus, $(v_k, y_k) \in \partial f(x_k,u_k)$ for each $k \in \N$ while,
	along a subsequence (without relabeling), $\norm{y_k}^{-1}(v_k, y_k) \to (0,\hat y)$
	for some $\hat y \in \R^m$ with $\hat y\neq 0$.
	As \cite[Proposition 3.2(a)]{LevyPoliquinRockafellar2000} yields 
	the convergence $f(x_k,u_k) \to f(\bar x, \bar u)$,
	we get $(0,\hat y) \in \partial^\infty f(\bar x, \bar u)$, violating \eqref{eq:basic_CQ}
	and, thus, \eqref{eq:basic_CQ_composite}.
	\end{proof}

\begin{remark}
	We have noticed that the proofs of \cite[Proposition~2.2]{LevyPoliquinRockafellar2000}
	and \cite[Proposition~13.32]{RockafellarWets1998} as well as the more recent result \cite[Theorem~1.61]{Mordukhovich2024},
	all justifying continuous prox-regularity of $f$ in the composite form from \cref{ass:composite_setting},
	do not really prove subdifferential continuity of $f$ and just state that it can be done easily.
	However, the straightforward arguments (the authors likely had in mind) do not work for that purpose.
	On the one hand, this potential gap can be overcome by assuming continuity of $g$ relative to its domain,
	which hardly has any restrictive impact on applications.
	On the other hand, we are aware of efforts (being not yet published)
	to bypass the need for subdifferential continuity.	
\end{remark}

In order to study the effects of perturbations induced by parameters $(v,u)\in\R^n\times\R^m$, we define the mappings
$M\colon\R^n\times\R^m\tto\R^n$ and $M_1\colon\R^n\times\R^m\tto\R^n\times\R^m$ by
\begin{align*}
		M(v,u)
		& :=
		\{x \,|\, v \in \partial_x f(x,u)\},
		\\
		M_1(v,u)
		& :=
		\{(x,y)\,|\,(v,y) \in \partial f(x,u)\}.
\end{align*}
Clearly, up to a permutation of variables, 
the graphs of $M$ and $M_1$ agree 
with the graphs of $\partial_x f$ and $\partial f$, respectively.
Moreover, $M_1$ has the same graph as the multiplier mapping $Y$.
\cref{Pro:Basic_properties_from_assumptions} yields that, in the composite setting from \cref{ass:composite_setting},
the graphs of $M_1$ and $Y$ are closed, and $\dom Y \subset \gph M$ holds locally around $(\bar v, \bar u, \bar x)$.
Under the qualification condition \eqref{eq:basic_CQ_composite}, 
this strengthens into $\dom Y = \gph M$, and both sets are locally closed around this reference triplet.

\begin{remark}\label{Rem:full_stab}
	As hinted in \cite[Theorem~2.3]{LevyPoliquinRockafellar2000}
	and explicitly stated in \cite[Corollary~1.3]{BenkoRockafellar2023},
	the so-called \emph{full stability} of a local minimizer $\bar x\in\R^n$ of \eqref{eq:Pbar},
	see \cite[Definition~1.1]{LevyPoliquinRockafellar2000} for the definition,
	implies that the mapping $M$ has a single-valued Lipschitz continuous 
	localization around $((\bar v,\bar u),\bar x)$.
	Similarly, the novel primal-dual full stability from 
	\cite[Definition~1.4]{BenkoRockafellar2023}
	implies that the mapping $M_1$ 
	has a single-valued Lipschitz continuous localization around $((\bar v,\bar u),(\bar x,\bar y))$,
	where $\bar y \in Y(\bar v,\bar u,\bar x)$, see \cite[Corollary~1.6]{BenkoRockafellar2023}.
	Thus, full-stability and its primal-dual version imply the isolated calmness of $M$ and $M_1$
	at (even around) $((\bar v,\bar u),\bar x)$ and $((\bar v,\bar u),(\bar x,\bar y))$, respectively,
	see \cref{sec:set_valued_maps} as well.
		Note that in the composite setting from \cref{ass:composite_setting},
		these implications are valid even without explicitly assuming \eqref{eq:basic_CQ_composite}
		because it was shown in the proof \cite[Theorem~5.1]{MordukhovichSarabi2018}
		that full stability automatically yields \eqref{eq:basic_CQ_composite}
		(that part of the proof does not use additional requirements on $g$ imposed in \cite[Theorem~5.1]{MordukhovichSarabi2018}),
		and, due to \cite[Proposition~2.2]{LevyPoliquinRockafellar2000},
		continuous prox-regularity as well.
\end{remark}

The mapping $M_1$ can be used to find stationary points of \eqref{eq:Pbar}
by the semismooth* Newton method stated in \cref{alg:SSN}
applied to solve the generalized equation
\begin{equation}\label{eq:GE_via_M1}
	(0,0) \in M^{-1}_1(x,y).
\end{equation}
Theoretically, $M$ could be used similarly, 
but its domain and image space do not have the same dimension.
Nevertheless, in the composite setting from \cref{ass:composite_setting}, 
$M$ can be replaced by
	the mapping $M_2\colon\R^n\times\R^m\tto\R^n\times\R^m$,
	motivated by \cite[Section~5]{GfrererOutrata2021} and given by
	\[
		M_2(v,u):=\{(x,w)\,|\,v\in\nabla f_0(x)+F'(x)^\top\partial g(w),\,u=w-F(x)\},
	\]
	which decouples the nonlinearities of $F$ and $g$ to some extent.
	Note the simple relation between the graphs of $M$ and $M_2$,
	valid for all $(v,u,x)$ near $(\bar v, \bar u, \bar x)$ under \eqref{eq:basic_CQ_composite}, namely
\begin{align*}
	((v,u),x)\in\gph M	
	\quad \iff	\quad
	((v,u),(x,F(x)+u))\in\gph M_2,
\end{align*}
since \cref{Pro:Basic_properties_from_assumptions} yields the description
\begin{equation*}
	M(v,u)=\{x\,|\,v\in\nabla f_0(x)+F'(x)^\top\partial g(F(x)+u)\}
\end{equation*}
of $M$. 
Since the multiplier $y$ does not explicitly appear in the descriptions 
of the mappings $M$ and $M_2$, we will refer to them as multiplier-free,
while we will call $M_1$ and $Y$ multiplier-based, since their common graph is characterized by \eqref{eq:perturbed_stationarity_system}.
Clearly, \cref{alg:SSN} can also be applied to solve the generalized equation
\begin{equation}\label{eq:GE_via_M2}
	(0,0) \in M^{-1}_2(x,w)
\end{equation}
in order to find stationary points of \eqref{eq:Pbar}.

Unlike $M_1$, $M_2$ is equivalent to $M$ in terms of various Lipschitzian properties.
In the next subsection, we study the isolated calmness (on a neighborhood) of these mappings, 
crucial for the superlinear convergence of \cref{alg:SSN}, 
see \cref{The:ss*_superlinear_convergence}.

\subsection{Variational analysis of perturbation mappings}
\label{sec:variational_analysis_perturbation_maps}

In order to examine various stability properties of the
mapping $M$, $M_1$, $M_2$, and $Y$,
we consider the closely related generalized derivatives of these mappings.
While we also provide some results regarding other stability notions,
let us recall that our focus is on the isolated calmness and its robust version.

Essentially, this subsection consists of three parts.
First, in the composite setting from \cref{ass:composite_setting}, 
we consider the relations between the multiplier-free mappings $M$ and $M_2$.
We will show that generalized derivatives as well as stability properties of these mappings
are related in a very straightforward manner.
Afterwards, we move on to establish a close connection between
derivatives and stability properties 
of the multiplier-based mappings $M_1$ and $Y$ which have the same graph.
These relations are valid in general, but we will also show that,
in the composite setting from \cref{ass:composite_setting},
we can explicitly compute the generalized derivatives of these mappings,
which, on the one hand, is a big advantage.
On the other hand, the disadvantage of relying on the multiplier-based mappings is that
various stability properties of these mappings are more restrictive than the corresponding
properties of the multiplier-free ones.
This will be confirmed in the third part, where, 
based on \cref{lem:calculus_for_limiting_tangent_cone},
we bridge the multiplier-free setting with the multiplier-based one.

\paragraph{Multiplier-free mappings}

Let us begin with the analysis of the multiplier-free mappings 
in the composite setting from \cref{ass:composite_setting} 
in the presence of the qualification condition \eqref{eq:basic_CQ_composite}.
Due to
\begin{equation}\label{eq:gph_M2} 
	\gph M_2 = \{((v,u),(x,w)) \,|\, (((v,u),x),F(x) + u - w) \in \gph M \times \{0\} \},
\end{equation}
(valid locally), $\gph M_2 $ is closed locally around $((\bar v,\bar u),(\bar x,F(\bar x)+\bar u))$ 
since $\gph M$ is locally around $((\bar v,\bar u),\bar x)$ by \cref{Pro:Basic_properties_from_assumptions}.
Moreover,
we can apply the change-of-coordinates formulas from \cite[Exercise~6.7]{RockafellarWets1998}, 
\cite[Proposition~4.4]{GfrererOutrata2023}, and \cref{lem:simple_calculus_for_tangents} 
to easily find the following result.

\begin{lemma}\label{lem:tangents_and_normals_to_M_and_M2}
	Let \cref{ass:composite_setting}
	and the qualification condition \eqref{eq:basic_CQ_composite} hold.
	For each $((v, u), x) \in \gph M$ close enough to $((\bar v,\bar u),\bar x)$, 
	the following assertions are valid.
	\begin{enumerate}
		\item Fix $(d^v,d^u,d^x)\in\R^n\times\R^m\times\R^n$.
			Then, for each of the tangent cone operators $\mathcal T\in\{T,T^{\#},T^\textup{P}\}$, 
			the following statements are equivalent:
			\begin{enumerate}
				\item[(i)] $((d^v,d^u),d^x)\in\mathcal T_{\gph M}(( v, u), x)$,
				\item[(ii)] $((d^v,d^u),(d^x,F'( x)d^x+d^u))\in\mathcal T_{\gph M_2}(( v, u),( x,F( x)+ u))$.
			\end{enumerate}
		\item Fix $(\eta^v,\eta^u,\eta^x)\in\R^n\times\R^m\times\R^n$.
			Then, for each of the normal cone operators $\mathcal N\in\{\widehat N,N\}$,
			the following statements are equivalent:
			\begin{enumerate}
				\item[(i)] $((\eta^v,\eta^u),\eta^x)\in\mathcal N_{\gph M}(( v, u), x)$,
				\item[(ii)] $((\eta^v,\eta^u-\eta^y),(\eta^x-F'( x)^\top\eta^y,\eta^y))\in\mathcal N_{\gph M_2}(( v, u),( x,F( x)+ u))$
					holds for all $\eta^y\in\R^m$.
			\end{enumerate}
	\end{enumerate}
\end{lemma}

Note that we also have
\[
	\gph M = \{((v,u),x)\,|\,((v,u),(x,F(x)+u))\in\gph M_2\},
\]
i.e., $\gph M$ is a preimage of $\gph M_2$,
but the change-or-coordinates formulas cannot be applied in this
situation as the derivative of the appearing smooth mapping
never has full row rank.

\begin{corollary}\label{cor:Lipschitzian_properties_M_vs_M2}
	Let \cref{ass:composite_setting} and the qualification condition \eqref{eq:basic_CQ_composite} hold.
	For each $((v, u), x) \in \gph M$ close enough to $((\bar v,\bar u),\bar x)$,
	$M$ is isolatedly calm at (is isolated calm around,
	has the Aubin property at, has a single-valued Lipschitzian localization at)
	$((v, u), x)$
	if and only if
	$M_2$ is isolatedly calm at (is isolated calm around,
	has the Aubin property at, has a single-valued Lipschitzian localization at) 
	$((v, u),(x,F( x)+ u))$.
\end{corollary}
\begin{proof}
	For the statement about isolated calmness at the reference point, 
	one has to show
	\[
		DM(( v, u), x)(0,0)=\{0\}
		\;\Longleftrightarrow\;
		DM_2(( v, u),( x,F( x)+ u))(0,0)=\{(0,0)\},
	\]
	for the one about isolated calmness around the reference point,
	\[
		D^{\#}M(( v, u), x)(0,0)=\{0\}
		\;\Longleftrightarrow\;
		D^{\#}M_2(( v, u),( x,F( x)+ u))(0,0)=\{(0,0)\}
	\]
	has to be verified, 
	and the assertion about the single-valued Lipschitzian
	localization requires (among others) to prove
	\begin{equation}\label{eq:paratingent_equivalence}
		D^{\textup{P}}M(( v, u), x)(0,0)=\{0\}
		\;\Longleftrightarrow\;
		D^{\textup{P}}M_2(( v, u),( x,F( x)+ u))(0,0)=\{(0,0)\}.
	\end{equation}
	These equivalences, however, follow trivially since
	$((0,0),d^x)\in\mathcal T_{\gph M}(( v, u), x)$ equals
	$((0,0),(d^x,F'(x)d^x))\in\mathcal T_{\gph M_2}(( v, u),( x,F( x)+ u))$
	for each of the tangent cone operators $\mathcal T\in\{T,T^{\#},T^{\textup{P}}\}$ 
	due to \cref{lem:tangents_and_normals_to_M_and_M2}.
	
	The assertion about the Aubin property follows if we can can show
	\begin{equation}\label{eq:normal_equivalence}
		D^*M(( v, u), x)(0)=\{(0,0)\}
		\;\Longleftrightarrow\;
		D^*M_2(( v, u),( x,F( x)+ u))(0,0)=\{(0,0)\},
	\end{equation}
	but this is a direct consequence of \cref{lem:tangents_and_normals_to_M_and_M2} as well
	since it also provides the equivalence of
	$((\eta^v,\eta^u),0)\in N_{\gph M}(( v, u), x)$
	and, for all $\eta^y\in\R^m$,
	$((\eta^v,\eta^u-\eta^y),(-F'( x)^\top\eta^y,\eta^y))\in N_{\gph M_2}((v,u),(x,F(x)+u))$.
	
	The statement about the existence of a single-valued Lipschitzian localization
	follows combining \eqref{eq:paratingent_equivalence} and \eqref{eq:normal_equivalence}.
\end{proof}

\paragraph{Multiplier-based mappings}

We proceed with the multiplier-based mappings.
	Again, we consider the composite setting from \cref{ass:composite_setting},
	but we abstain from postulating validity of the qualification condition \eqref{eq:basic_CQ_composite} now.

\begin{lemma}\label{lem:tangents_to_gph_Y}
	Let \cref{ass:composite_setting} hold.
	For each $(( v, u),( x, y))\in\gph M_1$, the following assertions are valid.
	\begin{enumerate}
	\item Fix $(d^v,d^u,d^x,d^y)\in\R^n\times\R^m\times\R^n\times\R^m$.
	Then, for each of the tangent cone operators $\mathcal T\in\{T,T^{\#},T^\textup{P}\}$, 
	the following statements are equivalent:
	\begin{enumerate}
	\item[(i)] $((d^v,d^u),(d^x,d^y))
	\in \mathcal T_{\gph M_1}(( v, u),( x, y))$,
	\item[(ii)] $((d^v,d^u,d^x),d^y)
	\in \mathcal T_{\gph Y}(( v, u, x), y)$,
	\item[(iii)] $\nabla^2_{xx}L(x,y)d^x+F'(x)^\top d^y-d^v = 0$,
	$(F'(x)d^x+d^u,d^y)
	\in \mathcal T_{\gph(\partial g)}(F( x)+ u, y)$.
	\end{enumerate}
	\item  Fix $(\eta^v,\eta^u,\eta^x,\eta^y)\in\R^n\times\R^m\times\R^n\times\R^m$.
	Then, for each of the normal cone operators $\mathcal N\in\{\widehat N,N\}$,
	the following statements are equivalent:
	\begin{enumerate}
	\item[(i)] $((\eta^v,\eta^u),(\eta^x,\eta^y))
	\in \mathcal N_{\gph M_1}(( v, u),( x, y))$,
	\item[(ii)] $((\eta^v,\eta^u,\eta^x),\eta^y)
	\in \mathcal N_{\gph Y}(( v, u, x), y)$,
	\item[(iii)] $\nabla^2_{xx}L( x, y)\eta^v-F'( x)^\top \eta^u+\eta^x = 0$,
	$(\eta^u,F'( x)\eta^v+\eta^y)
	\in \mathcal N_{\gph(\partial g)}(F( x)+ u, y)$.
	\end{enumerate}
	\end{enumerate}
\end{lemma}
\begin{proof}
	The equivalencies between the statements (i) and (ii) are valid in general as we have $\gph M_1=\gph Y$.
	Equivalence to statement (iii) can be shown with the
	aforementioned change-of-coordinates formulas based on the representation
	\[
		\gph M_1
		=
		\{((v,u),(x,y))\,|\,
		(\nabla_xL(x,y)-v,F(x)+u,y)\in\{0\}\times\gph(\partial g)\}
	\]
	which follows from \eqref{eq:perturbed_stationarity_system}.
\end{proof}

Given the reference triplet $(\bar v,\bar u,\bar x)$, 
let us also fix $\bar y\in Y(\bar v,\bar u,\bar x)$.
It is clear that each of the primal derivative criteria 
\eqref{eq:LevyRockafellar_criterion}, \eqref{eq:GfrererOutrata_criterion},
or \eqref{eq:Strict_graph_der_criterion}, applied to the mapping $M_1$, which can be written as
\begin{equation}\label{eq:primal_crit_M_1}
	((0,0),(d^x,d^y))
	\in
	\mathcal T_{\gph M_1}((\bar v,\bar u),(\bar x,\bar y))
	\quad\Longrightarrow\quad
	(d^x,d^y) = (0,0)
\end{equation}
for the corresponding tangent cone operator $\mathcal T\in\{T,T^{\#},T^\textup{P}\}$,
can be decomposed into the two conditions
\begin{subequations}
	\begin{align}\label{eq:primal_crit_M}
		((0,0),(d^x,d^y))
		\in
		\mathcal T_{\gph M_1}((\bar v,\bar u),(\bar x,\bar y))
		&\quad\Longrightarrow\quad
		d^x = 0,
		\\\label{eq:primal_crit_Y}
		((0,0),(0,d^y))
		\in
		\mathcal T_{\gph M_1}((\bar v,\bar u),(\bar x,\bar y))
		&\quad\Longrightarrow\quad
		d^y = 0.
	\end{align}
\end{subequations}
Due to \cref{lem:tangents_to_gph_Y}, \eqref{eq:primal_crit_Y} is equivalent to
\[
	((0,0,0),d^y)
	\in
	\mathcal T_{\gph Y}((\bar v,\bar u,\bar x),\bar y)
	\quad\Longrightarrow\quad
	d^y = 0
\] 
and, thus, precisely coincides with the primal derivative criterion 
\eqref{eq:LevyRockafellar_criterion}, \eqref{eq:GfrererOutrata_criterion},
or \eqref{eq:Strict_graph_der_criterion} applied to the mapping $Y$.
An interpretation of \eqref{eq:primal_crit_M} will be provided in \cref{prop:IC_M}.

For our purposes, only the relations between the criteria 
\eqref{eq:LevyRockafellar_criterion} and \eqref{eq:GfrererOutrata_criterion}
for $M_1$ and $Y$ are important, showing that the isolated calmness of $M_1$ at (around) a point
implies the isolated calmness of $Y$ at (around) the same point.
The strict graphical derivative criterion \eqref{eq:Strict_graph_der_criterion} for $M_1$ and $Y$
plays a crucial rule in the recent paper \cite{BenkoRockafellar2023} 
dedicated to so-called \emph{primal-dual full stability}.

The isolated calmness of $Y$ at $((\bar v,\bar u,\bar x),\bar y)$, 
corresponding to \eqref{eq:primal_crit_Y} for $\mathcal T:=T$, deserves additional attention.
This condition particularly implies that the multiplier $\bar y$ associated with the
stationary point $\bar x$ of \hyperref[eq:P(u,v)]{(P$(\bar v,\bar u)$)} is uniquely determined
as $Y$ is convex-valued by \cref{Pro:Basic_properties_from_assumptions}\,\ref{item:Y_convex_valued}.
Since the isolated calmness of $M_1$ at $((\bar v,\bar u),(\bar x,\bar y))$, encoded by \eqref{eq:primal_crit_M_1}, 
implies \eqref{eq:primal_crit_Y}, it also implies uniqueness of the multiplier.

More can be said about the isolated calmness of $Y$ if we explore the following explicit form
of the the conditions
\eqref{eq:primal_crit_M_1}, \eqref{eq:primal_crit_M}, and \eqref{eq:primal_crit_Y}:
\begin{subequations}\label{seq:3_explicit_cond}
	\begin{align}\label{eq:IC_M1}
		\left.
		\begin{aligned}
		&\nabla^2_{xx}L(\bar x,\bar y)d^x+F'(\bar x)^\top d^y =0,
		\\
		&(F'(\bar x)d^x,d^y) \in \mathcal T_{\gph(\partial g)}(F(\bar x),\bar y)
		\end{aligned}
		\right\}
		&\quad\Longrightarrow\quad
		d^x=0,\,d^y=0,
		\\\label{eq:criticality}
		\left.
		\begin{aligned}
		&\nabla^2_{xx}L(\bar x,\bar y)d^x+F'(\bar x)^\top d^y =0,
		\\
		&(F'(\bar x)d^x,d^y) \in \mathcal T_{\gph(\partial g)}(F(\bar x),\bar y)
		\end{aligned}
		\right\}
		&\quad\Longrightarrow\quad 
		d^x=0,
		\\\label{eq:local_uniqueness_multiplier}
		F'(\bar x)^\top d^y=0,\,
		(0,d^y) \in \mathcal T_{\gph(\partial g)}(F(\bar x),\bar y)
		&\quad\Longrightarrow\quad
		d^y=0.
	\end{align}
\end{subequations}
For brevity of notation, let $\widehat Y(\bar x):=Y(\bar v,\bar u,\bar x)$
denote the set of multipliers associated with the pair $(\bar v,\bar u)$
of reference parameters and the fixed stationary point $\bar x$ of our interest.
Consider the mapping $Y_{\bar x}\colon \R^n \times \R^m \tto \R^m$ given by
$Y_{\bar x}(v,u):= Y(v,u,\bar x)$.
Based on the representation
\begin{equation*}
	Y_{\bar x}(v,u)
	=
	\{y\,|\,\nabla f_0(\bar x)+F'(\bar x)^\top y-v=0,\,(F(\bar x)+u,y)\in\gph(\partial g)\},
\end{equation*}
and \cite[Exercise~6.7]{RockafellarWets1998}, one can easily show that
\[
	((d^v,d^u),d^y)\in T_{\gph Y_{\bar x}}((\bar v,\bar u),\bar y)
	\;\Longleftrightarrow\;
	F'(\bar x)^\top d^y-d^v=0,\,(d^u,d^y)\in T_{\gph(\partial g)}(F(\bar x),\bar y),
\]
see \cref{lem:tangents_to_gph_Y} as well.
Hence, by the Levy--Rockafellar criterion, $Y_{\bar x}$ is isolatedly calm
at $((\bar v,\bar u),\bar y)$ if and only if condition \eqref{eq:local_uniqueness_multiplier}
holds for $\mathcal T:=T$. The latter, however, is already equivalent to the isolated calmness of $Y$ at
$((\bar v,\bar u,\bar x),\bar y)$.
All these conditions imply that $\widehat Y(\bar x)$ is a singleton.
In turn, if the multiplier $\bar y$ is uniquely determined and if $Y_{\bar x}$ is calm at
$((\bar v,\bar u),\bar y)$, then there are a neighborhood $V\subset\R^m$ of $\bar y$ 
and a constant $\kappa>0$ such that
\[
	\norm{y-\bar y}
	\leq
	\kappa(\norm{\nabla_xL(\bar x,y)}+\dist(F(\bar x),(\partial g)^{-1}(y)))
\]
holds for all $y\in V$, i.e., $Y_{\bar x}$ is isolatedly calm at $((\bar v,\bar u),\bar y)$.
Related observations have been formulated in 
\cite[Corollary~3.14]{DeMarchiMehlitz2023},
\cite[Theorem~4.5]{HangMordukhovichSarabi2020},
\cite[Theorem~8.1]{MohammadiMordukhovichSarabi2022}, and
\cite[Theorem~4.1]{MordukhovichSarabi2019}.
We summarize the above observations in the subsequently stated corollary.

\begin{corollary}\label{cor:uniqueness_result}
	Let \cref{ass:composite_setting} hold.
	For fixed $\bar y\in\widehat Y(\bar x)$, the following statements are equivalent.
	\begin{enumerate}
		\item The mapping $Y$ is isolatedly calm at $((\bar v,\bar u,\bar x),\bar y)$.
		\item The mapping $Y_{\bar x}$ is isolatedly calm at $((\bar v,\bar u),\bar y)$.
		\item We have $\widehat Y(\bar x)=\{\bar y\}$,
			and $Y_{\bar x}$ is calm at $((\bar v,\bar u),\bar y)$.
		\item Condition \eqref{eq:primal_crit_Y}, or explicitly \eqref{eq:local_uniqueness_multiplier}, 
			holds for $\mathcal T:=T$.
	\end{enumerate}
\end{corollary}

\paragraph{Connecting the multiplier-free and multiplier-based mappings}

Finally, we aim to connect the isolated calmness property of the multiplier-free mapping $M$
with the isolated calmness of the multiplier-based mappings $M_1$ and $Y$.
This connection is valid in general, 
so we proceed without assuming the composite structure from \cref{ass:composite_setting},
while noting that with this special structure, the connection extends also to the mapping
$M_2$ on the basis of \cref{cor:Lipschitzian_properties_M_vs_M2}.
Interestingly, we only need to apply \cref{lem:calculus_for_limiting_tangent_cone}
to the mapping $Y$.

\begin{lemma}\label{lem:graphical_derivative_M_vs_Y}
	Let $f$ be chosen such that the properties~\ref{item:Y_convex_valued} and \ref{item:Y_gph_closed} 
	from \cref{Pro:Basic_properties_from_assumptions} hold.
	Then, for $\mathcal T := T$ and $\mathcal D := D$ 
	($\mathcal T := T^{\#}$ and $\mathcal D := D^{\#}$), we have
	\begin{equation}\label{eq:graphical_derivative_M_vs_Y}
		((d^v,d^u),d^x) \in \mathcal T_{\gph M}((\bar v,\bar u),\bar x)
		\quad \Longleftarrow \quad
		(d^v,d^u,d^x) \in \bigcup\limits_{\bar y\in \widehat Y(\bar x)}
		\dom \mathcal D Y((\bar v,\bar u,\bar x),\bar y).
	\end{equation}
	Moreover, equivalence holds if also the properties 
	\hyperref[item:Y_dom_closed]{(c)} and \hyperref[item:Y_loc_bounded]{(d)} 
	from \cref{Pro:Basic_properties_from_assumptions} are satisfied 
	and if $Y$ is (locally uniformly) inner calm* in the fuzzy sense 
	w.r.t.\ $\dom Y$ at $(\bar v,\bar u,\bar x)$ in direction $(d^v,d^u,d^x)$.
\end{lemma}
\begin{proof}
Thanks to property~\ref{item:Y_convex_valued}, 
we have $\mathcal T_{\dom Y}(\bar v,\bar u,\bar x) \subset \mathcal T_{\gph M}((\bar v,\bar u),\bar x)$
while property~\ref{item:Y_gph_closed} enables us to use \cref{lem:calculus_for_limiting_tangent_cone}, 
showing the first statement.
For the second statement, property~\hyperref[item:Y_dom_closed]{(c)} turns the above inclusion into equality 
while property~\hyperref[item:Y_loc_bounded]{(d)} together with the (locally uniform) fuzzy inner calmness* of $Y$ 
guarantees also validity of the converse implication in \cref{lem:calculus_for_limiting_tangent_cone}.
\end{proof}

Now, we are in position to provide the promised interpretation of condition \eqref{eq:primal_crit_M},
which, in combination with a suitable inner calmness* assumption,
provides a characterization of the isolated calmness of $M$ at/around $((\bar v,\bar u),\bar x)$.

\begin{theorem}\label{prop:IC_M}
	Let $f$ be chosen such that the properties~\ref{item:Y_convex_valued} and \ref{item:Y_gph_closed} 
	from \cref{Pro:Basic_properties_from_assumptions} hold
	(exemplary given if \cref{ass:composite_setting} is valid).
	If the mapping $M$ is isolatedly calm at (around) $((\bar v,\bar u),\bar x)$, then
	\eqref{eq:primal_crit_M} for $\mathcal T:=T$ ($\mathcal T:=T^{\#}$)
	is valid for each $\bar y\in\widehat Y(\bar x)$,
	and $Y$ is (locally uniformly) inner calm* in the fuzzy sense
	w.r.t.\ $\dom Y$ at $(\bar v,\bar u,\bar x)$ in each direction
	$(0,0,d^x)$ with $d^x \neq 0$.
	The converse implication holds if also the properties 
	\hyperref[item:Y_dom_closed]{(c)} and \hyperref[item:Y_loc_bounded]{(d)} 
	from \cref{Pro:Basic_properties_from_assumptions} are satisfied 
	(exemplary given if \cref{ass:composite_setting} and the qualification condition \eqref{eq:basic_CQ_composite} are valid).
\end{theorem}
\begin{proof}
	The second statement follows from 
	\cref{lem:tangents_to_gph_Y,lem:graphical_derivative_M_vs_Y}
	as well as the primal derivative criteria 
	\eqref{eq:LevyRockafellar_criterion} and \eqref{eq:GfrererOutrata_criterion}.
	
	In order to justify the first statement, note that the isolated calmness assumption
	is sufficient for the validity of the derivative criterion
	\eqref{eq:primal_crit_M} for $\mathcal T:=T$ or $\mathcal T:=T^{\#}$
	even if $\gph M$ is not locally closed.
	It, thus, remains to argue why the respective inner calmness* assumption is also necessary.
	This, however, follows easily since, 
	for $\mathcal T:=T$ ($\mathcal T:=T^{\#}$),
	\eqref{eq:primal_crit_M} implies that 
	$((0,0),d^x) \in \mathcal T_{\gph M}((\bar v,\bar u),\bar x)$
	only for $d^x = 0$.
	Thus, for each $d^x \neq 0$, we have 
	$(0,0,d^x) \notin \mathcal T_{\dom Y}(\bar v,\bar u,\bar x)$
	and (locally uniform) inner calmness* in the fuzzy sense
	w.r.t.\ $\dom Y$ of $Y$ at $(\bar v,\bar u,\bar x)$ in direction $(0,0,d^x)$ follows,
	see \cref{Rem:uniform_ic*_automatic}.
\end{proof}

\begin{remark}\label{Rem:On_fic*}
Interestingly, the above \cref{prop:IC_M} shows that the inner calmness* assumptions,
needed to make the explicit conditions \eqref{eq:primal_crit_M} 
sufficient for the isolated calmness (on a neighborhood) of $M$,
are inherently satisfied under that isolated calmness.
This is comforting since we know that we do not add anything superficial.
Nevertheless, if we aim to show validity of the isolated calmness of $M$ (on a neighborhood)
using \eqref{eq:primal_crit_M}, we still need to find a way how to verify
these inner calmness* assumptions, which is not an easy task.
Structural properties of the problem data can provide some help.
As shown in \cite[Theorem~3.9(ii)]{Benko2021}, for NLPs,
(the basic, not necessarily uniform) inner calmness* in the fuzzy sense of $Y$ always holds,
and \eqref{eq:primal_crit_M} for $\mathcal T:=T$, thus, fully characterizes 
the isolated calmness of $M$ at $((\bar v,\bar u),\bar x)$.
Moreover, in the presence of \cref{ass:composite_setting}, 
taking into account that $(w,y)\in\gph(\partial g)$ is equivalent to 
$((w,g(w)),(y,-1))\in\gph N_{\epi g}$,
\cite[Theorem~3.9(ii)]{Benko2021} can be easily extended to the case where
$\epi g$ is a convex polyhedral set.
	The derivation of further sufficient criteria 
	for these inner calmness*
	conditions is an important subject of future research.
	Exemplary, we conjecture that so-called \emph{strict complementarity},
	demanding the existence of $\bar y\in Y(\bar x)$ belonging to the relative interior of $\partial g(F(\bar x))$
	and
	utilized recently e.g.\ in 
	\cite[Theorem~4.1]{HangJungSarabi2024},
	\cite[Theorems~3.10 and~3.12]{HangSarabi2024}, or
	\cite[Theorem~5.10]{MordukhovichSarabi2019}, yields fuzzy inner calmness*,
	but a detailed study of this supposition is beyond the capacity of this paper.
\end{remark}

We conclude this section by restating the equivalence between \eqref{eq:primal_crit_M_1} 
and the combination of \eqref{eq:primal_crit_M} and \eqref{eq:primal_crit_Y} 
for $\mathcal T:=T$ ($\mathcal T:=T^{\#}$) in terms of isolated calmness at (around) a point.
In the presence of properties~\ref{item:Y_convex_valued},~\ref{item:Y_gph_closed}, and~\hyperref[item:Y_loc_bounded]{(d)} 
from \cref{Pro:Basic_properties_from_assumptions},
\cref{lem:sufficient_conditions_IC*} yields that
\eqref{eq:primal_crit_Y} for $\mathcal T:=T$ ($\mathcal T:=T^{\#}$) 
is a sufficient condition for (locally uniform) fuzzy inner calmness* of $Y$.
Since \eqref{eq:primal_crit_M_1} implies \eqref{eq:primal_crit_Y}, 
it, too, provides such a sufficient condition.
	Thus, we obtain the following result as a corollary of \cref{prop:IC_M},
	noting that $\gph M_1$ is always closed whenever property~\ref{item:Y_gph_closed}
	of \cref{Pro:Basic_properties_from_assumptions} holds.

\begin{corollary}\label{prop:IC_M1_M_Y}
	Let $f$ be chosen such that the properties~\ref{item:Y_convex_valued} and \ref{item:Y_gph_closed} 
	from \cref{Pro:Basic_properties_from_assumptions} hold
	(exemplary given if \cref{ass:composite_setting} is valid).
	If the mapping $M$ is isolatedly calm at (around) $((\bar v,\bar u),\bar x)$ while
	the mapping $Y$ is isolatedly calm at (around) $((\bar v,\bar u,\bar x),\bar y)$,
	then the mapping $M_1$ is isolatedly calm at (around) $((\bar v,\bar u),(\bar x,\bar y))$.
	The converse implication holds if also the properties 
	\hyperref[item:Y_dom_closed]{(c)} and \hyperref[item:Y_loc_bounded]{(d)} 
	from \cref{Pro:Basic_properties_from_assumptions} are satisfied 
	(exemplary given if \cref{ass:composite_setting} and the qualification condition \eqref{eq:basic_CQ_composite} are valid).
\end{corollary}

\section{Critical multipliers and local fast convergence of Newton-type methods}
\label{sec:critical_multipliers}

Recall that generalized nonlinear programming corresponds to problem \eqref{eq:Pbar}
under the composite structure from \cref{ass:composite_setting}, which, in unperturbed form, reads
\begin{equation}\label{eq:CP}\tag{CP}
	\min\limits_x\quad f_0(x) + g(F(x))\quad\text{s.t.}\quad x\in\R^n.
\end{equation}
	Let us fix a reference triplet $(\bar v,\bar u,\bar x)\in\R^n\times\R^m\times\R^n$ 
	with $(\bar v,\bar u):=(0,0)$ and $\bar x$ being a stationary point of \eqref{eq:CP}
	as we are mainly interested in the properties of the unperturbed problem.
	Throughout the section,
	we abstain from denoting the reference parameters by $(\bar v,\bar u)$ for clarity of notation. 
	
To start the discussion about critical multipliers, 
let us consider the equality-constrained optimization problem
\begin{equation}\label{eq:P=}
	\min\limits_x\quad f_0(x)\quad\text{s.t.}\quad F(x)=0,
\end{equation}
which is a special case of \eqref{eq:CP} for $g:=\delta_{\{0\}}$.
Let us fix $\bar y\in \widehat Y(\bar x)$.
Following \cite{Izmailov2005}, $\bar y$ is said to be a
\emph{critical multiplier} whenever there exist a nonvanishing $d^x\in\R^n$
and some $d^y\in\R^m$ such that
\[
    \nabla^2_{xx} L(\bar x, \bar y) d^x + F'(\bar x)^\top d^y=0,\,
    F'(\bar x)d^x=0.
\]
Observe that the normal cone mapping associated with the set $\{0\}$ possesses
the graph $\{0\}\times\R^m$, so that we can rewrite the above in seemingly more
complicated form as
\[
	\nabla^2_{xx}L(\bar x,\bar y) d^x + F'(\bar x)^\top d^y=0,\,
	(F'(\bar x)d^x,d^y)\in T_{\gph N_{\{0\}}}(F(\bar x),\bar y),
\]
or, by means of the definition of the graphical derivative, as
\begin{equation}\label{eq:criticality_P=}
	\nabla^2_{xx}L(\bar x,\bar y) d^x + F'(\bar x)^\top d^y=0,\,
	d^y \in D N_{\{0\}}(F(\bar x),\bar y)(F'(\bar x)d^x).
\end{equation}

Let us now switch over to the slightly more general constrained optimization problem
\begin{equation}\label{eq:PD}
	\min\limits_x\quad f_0(x)\quad\text{s.t.}\quad F(x)\in D
\end{equation}
for some closed convex set $D\subset\R^m$.
Clearly, \eqref{eq:P=} is a special case of \eqref{eq:PD} for $D:=\{0\}$,
and \eqref{eq:PD} is a special case of \eqref{eq:CP} for $g:=\delta_D$.
Let us fix $\bar y\in \widehat Y(\bar x)$.
Motivated by \eqref{eq:criticality_P=}, we refer to $\bar y$ as critical
whenever there are a nonvanishing $d^x\in\R^n$ and some $d^y\in\R^m$ such that
\begin{equation}\label{eq:criticality_PD}
	\nabla^2_{xx}L(\bar x,\bar y) d^x + F'(\bar x)^\top d^y=0,\,
	d^y \in D N_{D}(F(\bar x),\bar y)(F'(\bar x)d^x),
\end{equation}
see \cite[Definition~3.1]{MordukhovichSarabi2019} as well.
In the special case where $D:=\R^m_-$ is the nonpositive orthant,
a simple calculation shows that 
the condition $d^y\in DN_D(F(\bar x),\bar y)(F'(\bar x)d^x)$
corresponds to demanding
\begin{equation}\label{eq:graphical_derivative_R_-}
	\begin{aligned}
		&F_i(\bar x)<0,\,\bar y_i=0&
		&\Longrightarrow&
		&d^y_i=0,&
		\\
		&F_i(\bar x)=0,\,\bar y_i>0&
		&\Longrightarrow&
		&F_i'(\bar x)d^x=0,&
		\\
		&F_i(\bar x)=0,\,\bar y_i=0&
		&\Longrightarrow&
		&F'_i(\bar x)d^x\leq 0,\,d^y_i\geq 0,\,F'_i(\bar x)d^x\,d^y_i=0&
	\end{aligned}
\end{equation}
for all $i\in\{1,\ldots,m\}$ in the given situation.
Particularly, we recover the definition of criticality from 
\cite[Definition~2]{IzmailovSolodov2012} which addresses 
NLPs with inequality (and equality) constraints.

Observe that \eqref{eq:criticality_PD} is the same as
\[
	\nabla^2_{xx}L(\bar x,\bar y) d^x + F'(\bar x)^\top d^y=0,\,
	d^y \in D(\partial\delta_D)(F(\bar x),\bar y)(F'(\bar x)d^x).
\]
This motivates the following definition of (non-)criticality which addresses
our model problem \eqref{eq:CP}.

\begin{definition}\label{def:criticality}
	Fix some multiplier $\bar y\in \widehat Y(\bar x)$.
	Then $\bar y$ is called a \emph{noncritical} multiplier of 
	\eqref{eq:CP} for $\bar x$ if 
	\begin{equation}\label{eq:non_criticality}
		\nabla^2_{xx}L(\bar x,\bar y)d^x + F'(\bar x)^\top d^y=0,\,
		d^y\in D(\partial g)(F(\bar x),\bar y)(F'(\bar x)d^x)
		\quad\Longrightarrow\quad
		d^x=0
	\end{equation}
	is valid. 
	Otherwise, we refer to $\bar y$ as \emph{critical} for $\bar x$.
\end{definition}

	We note that \eqref{eq:non_criticality} is equivalent to \eqref{eq:criticality}
	for $\mathcal T:=T$.
	The above definition of noncriticality can also be found in \cite[Definition~3.1]{MordukhovichSarabi2018}
	where it has been stated for a particular class of convex, piecewise affine functions $g$.
	Let us also mention that, in \cite{MordukhovichSarabi2018}, the slightly more
	general setting of a so-called perturbed \emph{variational system} of the
	form
	\begin{equation}\label{eq:perturbed_variational_system}
		J(x) + F'(x)^\top y = v,
		\qquad
		y\in\partial g(F(x)+u) 
	\end{equation}
	for a continuously differentiable mapping $J\colon\R^n\to\R^n$ has been considered.
	Clearly, \eqref{eq:perturbed_variational_system} 
	is closely related to the variational inclusion $-J(x)\in\partial (g\circ F)(x)$.
	The perturbed stationarity system \eqref{eq:perturbed_stationarity_system}
	that we are investigating here results from 
	\eqref{eq:perturbed_variational_system} by choosing $J:=\nabla f_0$,
	and one can easily check that our modeling approach suggested in
	\cref{sec:model} can be easily adapted to cover perturbed variational systems
	of type \eqref{eq:perturbed_variational_system}, 
	simply by continuous differentiability of $J$ and twice
	continuous differentiability of $f_0$, which guarantee that these mappings
	do not cause any difficulties in the variational calculus.
	However, as our main motivation is induced by the setting of composite optimization,
	we stick to the slightly more restrictive setting 
	discussed in \cref{sec:perturbed_optimization_problems}.

Let us fix some multiplier $\bar y\in \widehat Y(\bar x)$.
Assume for a moment that
\eqref{eq:basic_CQ_composite} holds and that
$g$ is also so-called \emph{twice epi-differentiable} at $F(\bar x)$ for $\bar y$,
see \cite[Definition~13.6]{RockafellarWets1998} for a precise definition.
Typical examples for proper, lower semicontinuous, convex functions $g$ 
satisfying the latter requirement
are the piecewise affine functions considered in \cite{MordukhovichSarabi2018} 
(comprising indicator functions of convex polyhedral sets)
or the indicator function of the second-order cone, 
see \cite[Theorem~3.1]{HangMordukhovichSarabi2020}.
Since $g$ is continuously prox-regular
as mentioned in the proof of \cref{Pro:Basic_properties_from_assumptions},
we have
\begin{equation*}
	D(\partial g)(F(\bar x),\bar y)(F'(\bar x)d^x)
	=
	\frac12\partial(\d^2 g(F(\bar x),\bar y))(F'(\bar x)d^x)
\end{equation*}
from \cite[Theorem~13.40]{RockafellarWets1998}, so that 
\[
	0
	\in \nabla^2_{xx}L(\bar x,\bar y)d^x
	+
	F'(\bar x)^\top
	D(\partial g)(F(\bar x),\bar y)(F'(\bar x)d^x)
\]
means, equivalently, that $d^x$ is a stationary point of
\begin{equation}\label{eq:stationarity_subproblem}
	\min\limits_s\quad
	\frac12\nabla^2_{xx}L(\bar x,\bar y)[s,s]
	+
	\frac12\d^2 g(F(\bar x),\bar y)(F'(\bar x)s)
	\quad\text{s.t.}\quad s\in\R^n
\end{equation}
provided the chain rule applies (which, particularly, is the case
when $\d^2 g(F(\bar x),\bar y)$ is a piecewise affine function).
Hence, in some situations, noncriticality of $\bar y$ means that
the all-zero vector is the uniquely determined stationary point 
of \eqref{eq:stationarity_subproblem}.
Let us note that this might be weaker than demanding that the
all-zero vector is the uniquely determined global minimizer
of \eqref{eq:stationarity_subproblem} which corresponds to
a second-order sufficient optimality condition for 
\eqref{eq:CP}, see \cite[Remark~3.15]{DeMarchiMehlitz2023}
(even for $g$ being merely lower semicontinuous).
Note that \cite[Section~3.3]{DeMarchiMehlitz2023} also shows that extending 
\cref{def:criticality} to situations where $g$ is not assumed to be convex
is also reasonable. 
Due to \cite[Theorem~5.2]{BenkoMehlitz2023}, already
\[
	\inf\limits_{s\in\R^n}\sup\limits_{y\in \widehat{Y}(\bar x)}
	\left(\nabla^2_{xx}L(\bar x,y)[s,s]
	+
	\d^2 g(F(\bar x),y)(F'(\bar x)s)\right)
	>
	0
\]
serves as a (weak) second-order sufficient optimality condition for \eqref{eq:CP},
and the latter, notably, does not rule out the existence of a critical multiplier.
Let us also mention that, in order to obtain second-order sufficient conditions
of reasonable strength, one would typically restrict the variable $s$ above to a suitable
critical cone, see \cite[Section~5]{BenkoMehlitz2023} for details.

Based on our analysis in \cref{sec:variational_analysis_perturbation_maps},
we are in position to characterize noncriticality of multipliers
in terms of the isolated calmness of the perturbation mappings studied therein.
Let us start with the following result which interrelates 
noncriticality with the isolated calmness of the mapping $M_1$.

\begin{corollary}\label{cor:non_criticality_via_IC_M1}
	Fix some multiplier $\bar y\in \widehat Y(\bar x)$.
	Then the following statements are equivalent.
	\begin{enumerate}
		\item The mapping $M_1$ is isolatedly calm at $((0,0),(\bar x,\bar y))$.
		\item The multiplier $\bar y$ is noncritical for $\bar x$,
			and \eqref{eq:local_uniqueness_multiplier} holds for $\mathcal T:=T$.
	\end{enumerate}
\end{corollary}
\begin{proof}
	We have seen in \cref{sec:variational_analysis_perturbation_maps} that
	the Levy--Rockafellar criterion for the mapping $M_1$ at
	$((0,0),(\bar x,\bar y))$ reduces to \eqref{eq:IC_M1} for $\mathcal T:=T$, 
	and that the latter condition decouples into \eqref{eq:criticality} as well as
	\eqref{eq:local_uniqueness_multiplier} for $\mathcal T:=T$, respectively.	
	Hence, the assertion follows from the closedness of $\gph M_1$.
\end{proof}

As we already outlined in \cref{cor:uniqueness_result},
condition \eqref{eq:local_uniqueness_multiplier} for $\mathcal T:=T$ implies uniqueness of the
underlying multiplier. However, the theory on critical multipliers is
typically employed to study situations where a problem-tailored 
(strong) second-order sufficient condition is violated and the
multiplier associated with the stationary point under consideration is not uniquely
determined. Hence, using isolated calmness of $M_1$ as a sufficient
condition for noncriticality might be of limited practical use.
A result similar to \cref{cor:non_criticality_via_IC_M1} 
has been shown in \cite[Theorem~7.5]{MordukhovichSarabi2018}.

Much more interesting is the following result,
a consequence of \cref{prop:IC_M}, which interrelates
noncriticality with the isolated calmness of $M$.

\begin{corollary}\label{cor:non_criticality_via_IC_M}
	Consider the following statements.
	\begin{enumerate}
		\item\label{item:IC_M_gives_non_criticality} The mapping $M$ is isolatedly calm at $((0,0),\bar x)$. 
		\item\label{item:non_criticality_gives_IC_M} 
			Each multiplier in $\widehat Y(\bar x)$ is noncritical for $\bar x$
			and $Y$ is inner calm* in the fuzzy sense w.r.t.\
			$\dom Y$ at $(0,0,\bar x)$ in each direction
			$(0,0,d^x)$ with $d^x \neq 0$.
	\end{enumerate}
		Then~\ref{item:IC_M_gives_non_criticality} always implies~\ref{item:non_criticality_gives_IC_M}
		while the converse implication holds in the presence of \eqref{eq:basic_CQ_composite}.
\end{corollary}

In contrast to \cref{cor:non_criticality_via_IC_M1}, 
\cref{cor:non_criticality_via_IC_M} covers situations where the
set of multipliers is not a singleton and which, thus, is of
much higher practical interest. 
The comments about the inner calmness* assumption from \cref{Rem:On_fic*}
apply here as well.
Particularly, for a broad class of problems modeled with functions $g$ 
possessing a convex polyhedral epigraph,
which cover NLPs 
but also the setting investigated in \cite{MordukhovichSarabi2018},
the nonexistence of critical multipliers is fully characterized
by the isolated calmness of $M$ at $((0,0),\bar x)$ 
if the qualification condition \eqref{eq:basic_CQ_composite} holds.

\begin{remark}
	\cref{cor:non_criticality_via_IC_M} together with \cref{Rem:full_stab} 
	yield that full stability
	implies that all multipliers are noncritical.
	This has been shown, with some effort and for a particular class of problems, 
	e.g.\ in \cite[Theorem~5.1]{MordukhovichSarabi2018}.
	Using our approach, however, this conclusion follows trivially 
	and in the very general setting \eqref{eq:GNLP_setting}.
	This underlines the importance of noticing that the isolated calmness of $M$
	excludes critical multipliers, the observation around which this paper is built.
\end{remark}

The following example presents a simple situation where $M$ is isolatedly calm
while the set of multipliers is not a singleton and, hence, $M_1$ does not
enjoy the isolated calmness property.
 
\begin{example}\label{ex:M_IC_multiplier_not_unique}
	Let us consider \eqref{eq:CP} with $n:=m:=2$ and
	\[
		f_0(x):=\tfrac12x_1^2+\tfrac12(x_2+1)^2,
		\quad
		F(x):=\begin{pmatrix}x_1^3-x_2\\-x_2\end{pmatrix},
		\quad
		g:=\delta_{\R^2_-},
	\]
	i.e., an inequality-constrained nonlinear optimization problem.
	One can easily check that its uniquely determined global minimizer
	$\bar x:=(0,0)^\top$ satisfies \eqref{eq:basic_CQ_composite}
	which reduces to the standard Mangasarian--Fromovitz constraint qualification
	in the present situation.
	Furthermore, we find $\widehat Y(\bar x)=\conv\{\mathtt e_1,\mathtt e_2\}$.
	Picking $\bar y\in\widehat Y(\bar x)$ and arbitrary $d^x,d^y\in\R^2$, 
	we consider the conditions
	\[
		\begin{pmatrix}0\\0\end{pmatrix}
		=
		\nabla^2_{xx}L(\bar x,\bar y)d^x+F'(\bar x)^\top d^y
		=
		\begin{pmatrix}d^x_1\\d^x_2-d^y_1-d^y_2\end{pmatrix}
	\]
	and
	\begin{align*}
		d^y
		&\in
		D(\partial\delta_{\R^2_-})(F(\bar x),\bar y)(F'(\bar x)d^x)
		\\
		&=
		\left\{r\in\R^2\,\middle|\,
			\begin{aligned}
				&\bar y_i>0\quad\Longrightarrow\quad d^x_2=0\\
				&\bar y_i=0\quad\Longrightarrow\quad d^x_2\geq 0,\,r_i\geq 0,\,d^x_2r_i=0
			\end{aligned}
		\right\}.	
	\end{align*}
	For the evaluation of the graphical derivative, 
	we made use of \eqref{eq:graphical_derivative_R_-}.
	The first of these conditions yields $d^x_1=0$, 
	and as at least one component of $\bar y$ is positive, 
	the second condition implicitly requires $d^x_2=0$.
	Hence, each multiplier in $\widehat Y(\bar x)$ is noncritical for $\bar x$.
	Furthermore, as $\epi g$ is convex polyhedral, $Y$ is inner calm* in the fuzzy sense w.r.t.\
	its domain at $(0,0,\bar x)$ by
	\cite[Theorem~3.9(ii)]{Benko2021}.
	Now, \cref{cor:non_criticality_via_IC_M} shows 
	that $M$ is isolatedly calm at $((0,0),\bar x)$.
	We also note that, for each $\bar y\in\widehat Y(\bar x)$, $M_1$ is not isolatedly calm 
	at $((0,0),(\bar x,\bar y))$ as \eqref{eq:local_uniqueness_multiplier} does not hold for
	$\mathcal T:=T$, see \cref{cor:uniqueness_result,cor:non_criticality_via_IC_M1}.
\end{example}

In \cref{cor:non_criticality_via_IC_M}, 
we studied noncriticality of multipliers from a global kind of perspective
in terms of the multiplier-free mapping $M$, i.e., we formulated necessary
and sufficient criteria for noncriticality of \emph{all} multipliers in the multiplier set.
Similarly, \cref{cor:non_criticality_via_IC_M1} is a global result as the situation therein
implicitly requires that the multiplier set reduces to a singleton.
In \cite[Proposition~1]{IzmailovSolodov2012}, the authors characterize
noncriticality of \emph{some} multiplier in terms of an error bound condition,
and related results can be found in other papers,
see e.g.\ \cite[Theorem~4.1]{MordukhovichSarabi2018}, \cite[Theorem~5.6]{MordukhovichSarabi2019},
or \cite[Theorem~3.6, Proposition~3.8]{Sarabi2022}.
This offers a finer analysis, 
and these error bounds were used to establish local fast convergence 
of Newton-, SQP-, or multiplier-penalty-type methods in the literature
even in the case where critical multipliers exist.
In a forthcoming paper, we plan to adjust our approach to obtain this local kind of analysis
for critical multipliers for the general composite model \eqref{eq:CP},
and we will carve out some consequences of these findings for the local convergence of the
semismooth* Newton method from \cref{alg:SSN}.

Let us recall the generalized equations stated in \eqref{eq:GE_via_M1} and \eqref{eq:GE_via_M2}.
When applying \cref{alg:SSN} for their solution, 
it follows from \cref{The:ss*_superlinear_convergence} 
that isolated calmness of $M_1$ and $M_2$ around the point
$((0,0),(\bar x,\bar y))$ and $((0,0),(\bar x,F(\bar x))$ is a sufficient condition for
the method to be locally well-defined and superlinearly convergent when applied to
\eqref{eq:GE_via_M1} and \eqref{eq:GE_via_M2}, respectively, provided the mappings 
$M_1$ and $M_2$ are semismooth*
(which guarantees semismoothness* of $M_1^{-1}$ and $M_2^{-1}$).

Observing that isolated calmness on a neighborhood of a mapping can be checked in terms
of a criterion involving the limiting tangent cone to the graph, 
the following definition might be reasonable.
\begin{definition}\label{def:strong_criticality}
	Fix some multiplier $\bar y\in\widehat Y(\bar x)$.
	Then $\bar y$ is called a 
	\emph{strongly noncritical} multiplier of \eqref{eq:CP} for $\bar x$ if
	\begin{equation}\label{eq:strong_noncriticality}
		\nabla^2_{xx}L(\bar x,\bar y)d^x+F'(\bar x)^\top d^y=0,\,
		d^y\in D^{\#}(\partial g)(F(\bar x),\bar y)(F'(\bar x)d^x)
		\quad\Longrightarrow\quad
		d^x=0
	\end{equation}
	is valid. Otherwise, we refer to $\bar y$ as \emph{weakly critical} for $\bar x$.
\end{definition}

We note that, by definition, each strongly noncritical multiplier is also noncritical,
while each critical multiplier is also weakly critical. 
Moreover, \eqref{eq:strong_noncriticality} corresponds to \eqref{eq:criticality}
for $\mathcal T:=T^{\#}$.
Finally, due to \cref{lem:simple_calculus_for_tangents}, the concepts of
noncriticality (criticality) and strong noncriticality (weak criticality) coincide whenever
$\gph(\partial g)$ is the union of finitely many subspaces. This, exemplary, happens to be the
case whenever $g:=\delta_{\{0\}}$, i.e., in the case of equality-constrained programming,
where $\gph(\partial \delta_{\{0\}})=\{0\}\times\R^m$.
Taking a broader look at settings where $g$ is allowed to be nonconvex,
we also would like to mention $g:=\norm{\cdot}_0$ here, 
where $\norm{\cdot}_0$ counts the nonzero entries of the input vector, and
\[
	\gph(\partial\norm{\cdot}_0)
	=
	\{(u,y)\in\R^m\times \R^m\,|\,\forall i\in\{1,\ldots,m\}\colon\, u_iy_i=0\}.
\]
Let us note that $\norm{\cdot}_0$ has been addressed in \cite[Section~5.1]{DeMarchiMehlitz2023}
where the authors discuss sufficient conditions for local fast convergence of a
multiplier-penalty method applied to \eqref{eq:CP} with this particular function.

Let us mention that, for $g:=\delta_{\R^m_-}$, one can easily check that
\[
	d^y\in D^{\#}(\partial g)(F(\bar x),\bar y)(F'(\bar x)d^x)
\]
for directions $d^x\in\R^n$ and $d^y\in\R^m$ is equivalent to demanding
\[
	\begin{aligned}
		&F_i(\bar x)<0,\,\bar y_i=0&
		&\Longrightarrow&
		&d^y_i=0,&
		\\
		&F_i(\bar x)=0,\,\bar y_i>0&
		&\Longrightarrow&
		&F_i'(\bar x)d^x=0,&
		\\
		&F_i(\bar x)=0,\,\bar y_i=0&
		&\Longrightarrow&
		&F'_i(\bar x)d^x\,d^y_i=0&
	\end{aligned}
\]
for all $i\in\{1,\ldots,m\}$.
Taking \eqref{eq:graphical_derivative_R_-} into account,
\eqref{eq:strong_noncriticality} is strictly stronger than \eqref{eq:non_criticality}
in this comparatively simple setting of inequality-constrained optimization.

Below, we characterize strong noncriticality of multipliers in terms 
of isolated calmness properties of the mappings $M_1$ and $M$.
These corollaries can be obtained in similar fashion as 
\cref{cor:non_criticality_via_IC_M1,cor:non_criticality_via_IC_M}.

\begin{corollary}\label{cor:strong_non_criticality_via_IC_M1}
	Fix some multiplier $\bar y\in \widehat Y(\bar x)$.
	Then the following statements are equivalent.
	\begin{enumerate}
		\item The mapping $M_1$ is isolatedly calm around $((0,0),(\bar x,\bar y))$.
		\item The multiplier $\bar y$ is strongly noncritical for $\bar x$,
			and \eqref{eq:local_uniqueness_multiplier} holds for $\mathcal T:=T^{\#}$.
	\end{enumerate}
\end{corollary}

\begin{corollary}\label{cor:strong_non_criticality_via_IC_M}
	Consider the following statements.
	\begin{enumerate}
		\item\label{item:IC_M_gives_strong_non_criticality} The mapping $M$ is isolatedly calm around $((0,0),\bar x)$. 
		\item\label{item:strong_non_criticality_gives_IC_M}  
			Each multiplier in $\widehat Y(\bar x)$ is strongly noncritical for $\bar x$
			and $Y$ is locally uniformly inner calm* in the fuzzy sense w.r.t.\
			$\dom Y$ at $(0,0,\bar x)$ in each direction
			$(0,0,d^x)$ with $d^x \neq 0$.
	\end{enumerate}
	Then~\ref{item:IC_M_gives_strong_non_criticality} always implies~\ref{item:strong_non_criticality_gives_IC_M}
	while the converse implication holds in the presence of \eqref{eq:basic_CQ_composite}.
\end{corollary}

We note that the locally uniform fuzzy inner calmness* assumption in
statement~\ref{item:strong_non_criticality_gives_IC_M} of \cref{cor:strong_non_criticality_via_IC_M}
might be more restrictive compared to the fuzzy inner calmness* assumption in
statement~\ref{item:non_criticality_gives_IC_M} of \cref{cor:non_criticality_via_IC_M}.
Observing that $Y$ possesses convex images by convexity of $g$,
\cref{lem:sufficient_conditions_IC*,lem:tangents_to_gph_Y} show that
a sufficient condition for the required uniform fuzzy inner calmness* of $Y$
is given by \eqref{eq:local_uniqueness_multiplier} for $\mathcal T:=T^{\#}$.
The latter, however, implies that $\widehat Y(\bar x)$ is a singleton,
see \cref{cor:uniqueness_result}, 
and we anyway enter the restrictive setting of \cref{cor:strong_non_criticality_via_IC_M1}.
Let us also mention that we are not yet aware of a reasonably broad setting in which
the locally uniform inner calmness* would be automatically satisfied
(and \eqref{eq:local_uniqueness_multiplier} for $\mathcal T:=T^{\#}$ is not).
Nevertheless, as before, the combination of \cref{cor:strong_non_criticality_via_IC_M}
and \cref{Rem:full_stab} yields that full stability rules out the existence of weakly critical multipliers.

In the setting of equality-constrained optimization,
\eqref{eq:basic_CQ_composite} reduces to the well-known linear independence
constraint qualification (LICQ), and one can easily check, e.g.\ with
the aid of \cref{lem:simple_calculus_for_tangents} or by direct calculation, 
that \eqref{eq:local_uniqueness_multiplier} for $\mathcal T:=T^{\#}$
is implied by LICQ. Hence, in this setting and in the presence of
LICQ, isolated calmness of $M$ at and around $((0,0),\bar x)$
are the same and reduce to (strong) noncriticality of the associated
uniquely determined multiplier $\bar y\in \widehat Y(\bar x)$ which reads as
\[
	\nabla^2_{xx} L(\bar x,\bar y)d^x + F'(\bar x)^\top d^y=0,\,F'(\bar x)d^x=0
	\quad\Longrightarrow\quad
	d^x=0.
\]
Consulting \cref{cor:Lipschitzian_properties_M_vs_M2,cor:uniqueness_result}
as well as \cref{prop:IC_M1_M_Y}, we find that the mappings
$M_1$ and $M_2$ are also isolatedly calm around $((0,0),(\bar x,\bar y))$
and $((0,0),(\bar x,0))$, respectively.
Hence, \cref{The:ss*_superlinear_convergence} indicates that \cref{alg:SSN},
when applied to the generalized equations \eqref{eq:GE_via_M1} and \eqref{eq:GE_via_M2},
is likely to produce a sequence which converges superlinearly to
$(\bar x,\bar y)$ and $(\bar x,0)$, respectively, whenever it is initialized
close enough to these points.
Note that we have $M_1^{-1}(x,y)=(\nabla_xL(x,y),-F(x))$ in this particular
setting, i.e., $M_1^{-1}$ is a single-valued and continuously differentiable
mapping. Following \cite[comments after Algorithm~2]{GfrererOutrata2021},
the standard (local) Newton method applied to $M_1^{-1}$ corresponds to
an application of \cref{alg:SSN} to \eqref{eq:GE_via_M1},
and its local superlinear convergence under the noncriticality of $\bar y$
for $\bar x$ has already been worked out in 
	\cite[Section~4]{Izmailov2005} in the presence of LICQ,
	which particularly yields uniqueness of the multiplier $\bar y$.
	Hence, we recover this classical result.
	In \cite{IzmailovSolodov2012}, the authors show that mere noncriticality
	is enough to obtain superlinear convergence of a \emph{stabilized}
	version of the Newton method, called \emph{stabilized SQP}, 
	in equality-constrained optimization.
We refer the interested reader to the paper \cite{IzmailovSolodov2015}
where the authors present a nice overview of properties associated with critical and noncritical
multipliers in equality-constrained optimization.
Let us also note that, unlike \cref{alg:SSN}, the standard Newton method cannot be applied to $M_2^{-1}$ 
since the latter is not single-valued.

To conclude this section, we would like to point the reader to the fact that 
noncriticality of \emph{all} multipliers might not be enough to yield
local fast convergence of Newton-type methods in more general situations than
equality-constrained optimization, see e.g.\ \cite[Example~1]{IzmailovSolodov2012},
and this may even extend to the isolated calmness of $M_1$ and $M_2$ based on
\cref{cor:Lipschitzian_properties_M_vs_M2,cor:non_criticality_via_IC_M1,cor:non_criticality_via_IC_M}.
Indeed, \cref{The:ss*_superlinear_convergence} requires strong metric subregularity
around the reference point
of the mapping which appears in the generalized equation under consideration
(or, equivalently, 
isolated calmness around the reference point of the associated inverse mapping).
Exemplary, when addressing \eqref{eq:GE_via_M1} and \eqref{eq:GE_via_M2},
\cref{cor:strong_non_criticality_via_IC_M1,cor:strong_non_criticality_via_IC_M}
indicate that strong noncriticality of all multipliers is needed to ensure this,
and whenever $M_2$ is considered, the multiplier mapping, additionally, has to possess
the locally uniform fuzzy inner calmness* property.
As shown earlier, already in the context of inequality-constrained optimization,
strong noncriticality is a more restrictive concept than noncriticality in general.
In our future work, we aim to explore the concept of strong noncriticality even more,
and we also plan to derive more practicable sufficient conditions for the presence of
locally uniform fuzzy inner calmness* for the mapping $Y$ which,
in contrast to \eqref{eq:GfrererOutrata_criterion} from \cref{lem:sufficient_conditions_IC*},
do not already imply uniqueness of the multiplier.

Let us mention \cite[Example~1]{IzmailovSolodov2012} again, which is an interesting example
showing that the aforementioned stabilized SQP method for NLPs
may fail to work if inequality constraints are present, 
even if initialized close to a stationary point satisfying LICQ and the corresponding 
unique multiplier, because the subproblems may fail to possess a solution.
In that setting, however, $M_1$ is isolatedly calm around that stationary couple,
so that the semismooth* Newton method should work due to \cref{The:ss*_superlinear_convergence}.
This can indeed be seen from \cite[Example~5.13]{GfrererOutrata2021},
where an almost identical example is described,
and where it is shown that the semismooth* Newton method produces 
a superlinearly convergent sequence
(more precisely, a particular instance of the semismooth* Newton method 
with a specified approximation step is used therein),
while the Newton--Josephy method from \cite{Josephy1979} does not work 
as the appearing subproblems do not possess a solution.

\section{Concluding remarks}\label{sec:conclusions}

For a broad class of parametrized optimization problems \eqref{eq:P(u,v)} 
with the composite structure from \cref{ass:composite_setting},
we have characterized the isolated calmness at and around a point of the perturbation mapping $M$
in terms of an explicit condition and a calmness-type assumption.
We have also shown that the isolated calmness of $M$ around a point 
yields local superlinear convergence of a semismooth* Newton method 
when applied to an auxiliary mapping $M_2$ 
in order to find a stationary point of the unperturbed problem \eqref{eq:Pbar}.
Finally, we have derived a strong connection between the isolated calmness of $M$ at a point 
and nonexistence of critical multipliers. 
Particularly, these two conditions are equivalent for standard nonlinear programs
satisfying a qualification condition,
but also for other problems with inherent polyhedrality.

In a forthcoming paper, we plan to refine our approach to obtain an analogous connection
between noncriticality of a \emph{single} multiplier and a suitable isolated calmness assumption,
which should then extend the available characterizations of critical multipliers 
in terms of an error bound condition to the general composite model \eqref{eq:CP}.
Additionally, we will carve out some consequences of these findings 
for the local convergence of the semismooth* Newton method.

\subsection*{Acknowledgements}

The authors would like to thank Terry Rockafellar for valuable discussions, 
which inspired them to use the parametrized model problem 
from \cite{LevyPoliquinRockafellar2000}, 
the setting which fits perfectly with the ideas of the paper, 
particularly the crucial \cref{lem:graphical_derivative_M_vs_Y}.
Furthermore, the authors appreciate the comments of two anonymous reviewers
that helped to enhance the presentation of the obtained results.
The research of Mat\'u\v{s} Benko was supported by the Austrian Science Fund (FWF) 
under grant P32832-N
as well as by the infrastructure of the Institute of Computational Mathematics, 
Johannes Kepler University Linz, Austria.



\end{document}